\input amstex
\documentstyle{amsppt}
\magnification=\magstep1
\hsize=5in
\vsize=7.3in

\TagsOnRight  
\topmatter
\title Algebraic Coset Conformal field theories
\endtitle
\author Feng  Xu \endauthor

\address{Department of Mathematics, University of Oklahoma, 601 Elm Ave,
Room 423, Norman, OK 73019}
\endaddress
\email{xufeng\@ math.ou.edu}
\endemail
\abstract{            
All unitary Rational Conformal Field Theories (RCFT) are conjectured to
be related to unitary coset Conformal Field Theories, i.e., 
gauged Wess-Zumino-Witten
(WZW) models with compact gauge groups.
In this paper we  use
subfactor theory and ideas of algebraic quantum field theory 
to approach coset Conformal Field Theories.  Two  conjectures are
formulated  and their consequences are discussed.  Some results
are presented which prove the  conjectures in special cases.  In
particular, one of the  results states that a  class
of representations of coset $W_N$ ($N\geq 3$) algebras with
critical parameters are irreducible,
and under the natural compositions (Connes' fusion), they generate a
finite dimensional fusion ring whose structure constants
are completely determined, thus
proving a long-standing conjecture about the representations of these algebras
.}
\endabstract        
\thanks
%\thanks   
I'd  like
to thank the referee for helpful suggestions. This work is partially
supported by NSF grant DMS-9820935.
\endthanks 

\endtopmatter
\document

\heading \S1. Introduction \endheading
Conformal field theories (CFT) in two dimensions (cf. [MS]) have attracted
great
attention
among mathematicians in  recent years.  A large class of CFT known as
Rational CFT (RCFT) are
more amenable than general CFT and the classification of RCFT is an
outstanding open problem.   

Unitary coset CFT is a gauged Wess-Zumino-Witten (WZW) model with a
compact  Lie group $G$ as a gauge group, and 
$H$, a subgroup of $G$, is gauged (cf. [KS]).
It has been conjectured (cf. [MS],
[Witten])
that all unitary RCFT (except perhaps orbifolds, which
are relatively better understood and also similar
to the coset CFT, cf. Page 428 of [MS]) are related to unitary 
coset  CFT. 
In the literature there has been several different mathematical
approaches to CFT, see, for an example, [DL] for an approach by using
vertex operator algebras (cf. [B] or [FLM]).  
However, in the case of WZW model with
a compact gauge group $G$ , there is
a manifestly unitary formulation of these unitary CFT by using subfactor
theory (cf. [J], [W1] and [L1]) and ideas from algebraic
quantum field theory  (cf. [Haag]).
This formulation has various advantages besides producing the expected
results, see, for an example, [X1] and [X2] for some results on certain
rings which seem to be invisible in other approaches.  In view of
the importance of coset CFT among all RCFT, it is natural to use this
unitary formulation to study unitary coset CFT, and this is the
main purpose of this paper.    \par
To illustrate the new ideas in this paper we will focus on the case 
when $G$ is a simply connected semisimple compact Lie group of 
type $A$, i.e., $G=SU(N_1)\times SU(N_2)\times...\times SU(N_n)$. 
The ideas of this paper can be applied to all compact semisimple and
simply connected Lie 
groups and we plan
to consider them  in  separate publications. \par

We will describe in more details about the paper.  In subsection 2.1, We define
 conformal precosheaf  and their covariant 
representations  as in [GL].  
We then show how coset $H\subset G$  naturally
gives rise to an irreducible conformal precosheaf  (Prop. 2.2) 
together with a class of covariant
representations.  
As a first step we proved in Th. 2.3 an irreducibility result in
our setting by using
vertex operators (cf. [FZ] or [Kacv]). To do this we show in
Prop. 2.3 that certain
smeared vertex operators are
affiliated (cf. [Mv]) with some von Neumann algebras. Note that these 
smeared vertex operators are generally unbounded operators and
one must be careful with the formal manuplations of them. Prop. 2.3 is
proved by using a series of lemmas (lemma 1-6).
It can be seen from the proof that the vertex operator algebra for
coset $H\subset G$  as defined on P. 67 of  [Kacv] or 
\S5 of [FZ] can be thought as ``germs'' 
of our irreducible conformal precosheaf. \par    
Based on physicists' argument and well known examples, in 
subsection 2.3 two 
conjectures are formulated about these representations.  The first
conjecture states that the representations generate a finite dimensional
ring under certain compositions, in other words, 
the CFT is really ``rational".  The second
conjecture is a formula about the 
square root of the minimal index or the
statistical dimension (cf. 4.1) of these   
representations in terms of certain limits of characters, referred to
as ``Kac-Wakimoto" formula in [L5].  Both conjectures are highly nontrivial.
In particular, the second conjecture implies Kac-Wakimoto conjecture
in \S 2 of [KW]. \par
The rest of the paper is devoted to the proof of the two  conjectures
in special cases.

Section 3 is 
about proving finite index by using certain commuting squares which
play an important role in the type $II_1$ subfactor theory (cf. [Po], [We]
). However, we will consider factors of type $III_1$. First 
we introduced a notion of co-finiteness for a pair $H\subset G$. 
This notion is motivated by the conjectures of 2.3.  Prop. 3.1,
which follows from the commuting squares in lemma 3.1,  is
a key observation of this paper. It implies that if $H_1 \subset H_2 \subset 
G$ and $H_1 \subset G $ is cofinite, then 
$H_1 \subset H_2$ is cofinite. This leads to an infinite series of cofinite
pairs in Cor. 3.1. \par
Section 4 can be considered as an application of [X1] and [X2]  
given
the results of  3.1. We first recall some of the results of [X1] 
which are used in  4.2 to determine certain ring structures. These
results are  summarized in Th. 4.1. Prop. 4.2 is the key observation in
section 4. 
It allows one to study the representations of coset $H\subset G$ by using
the theory which has been developed for $G$ and $H$ in some cases. 
Then we prove Th. 4.2 which states that for a pair $H\subset G$
which is cofinite, if similar 
results  in Chap. V of [W2] hold for $H$,  then
Conj. 1 is also true. This result and Cor. 3.1 imply Conj. 1 for the
infinite series in Cor. 4.2. \par
In 4.3 we study the cosets corresponding to 
diagonal inclusions by using Th. 4.2 and
prove in Th. 4.3 that Conj. 2 is true for these cosets . We
also prove in Th. 4.3 certain irreducibility result. \par
The results in Th. 4.3 in the special case of coset $W_N$ algebras
with critical parameters 
have been long anticipated in both mathematics and
physics literature.  For example  these results are related to 
 the conjectures   in \S 3 and 4 of [FKW], which 
are  conjectured for what is know as ``W-algebras'', and they are closely
related to coset $W$-algebras in the vertex operator algebra sense (cf. 
[Watts]). \par 
%However the fusion products are not defined in the setting of
%[FKW] (cf. P. 2 of [FKW]).    \par
In 4.4 we present more examples by using Th. 4.2. The first
example is the coset $H\subset G$ where $H$ is the Cartan subgroup of $G$.
This coset is related to Parafermions (cf. [DL]).
We then consider a ``Maverick'' coset model considered in [DJ1] and 
determine the relevant ring structure. \par
All the groups considered in this paper are assumed to be 
connected compact Lie groups unless stated otherwise.  
As we already noted before, 
$G$ always denotes a simply connected semisimple compact Lie group of 
type $A$, i.e., $G=SU(N_1)\times SU(N_2)\times...\times SU(N_n)$. \par 
\heading \S2.  Coset CFT from Algebraic QFT point of view \endheading
\subheading {2.1 The irreducible conformal precosheaf  and its representations}
In this subsection we recall the basic properties enjoyed by the family of
the von Neumann algebras associated with a conformal quantum field theory
on $S^1$ (cf. [GL] and [FG]).               
By an {\it interval} in this paper 
 we shall always mean an open connected subset $I$
of $S^1$ such that $I$ and the interior $I' $ of its complement are
non-empty.  We shall denote by  ${\Cal I}$ the set of intervals in $S^1$.
We shall denote by  $PSL(2, {\bold R})$ the group of
  conformal transformations on the complex plane
that preserve the orientation and leave the unit circle $S^1$ globally
invariant.  Denote by ${\bold G}$
the universal covering group of $PSL(2, {\bold R})$.  Notice that  ${\bold G}$
is a simple Lie group and has a natural action on the  unit circle $S^1$. \par
Denote by  $R(\vartheta )$  the (lifting to ${\bold G}$ of the) rotation by
an angle $\vartheta $. This one-parameter subgroup of
${\bold G}$ will be referred to as rotation group (denoted by Rot) 
in the following.  
We may associate a 
one-parameter group with any interval $I$ in the following way.
Let $I_1$ be the upper
semi-circle, i.e. the interval
$\{e^{i\vartheta }, \vartheta \in (0, \pi )\}$.
 By using the Cayley transform
$C:S^1 \rightarrow {\bold R} \cup \{\infty \}$ given by
$z\rightarrow -i(z-1)(z +1)^{-1}$,
we may identify  $I_1$
with the positive real line ${\bold R}_+$. Then we consider
the one-parameter group $\Lambda _{I_1}(s)$ of diffeomorphisms of 
$S^1$  such that
$$
C\Lambda _{I_1} (s) C^{-1} x = e^s x \, ,
\quad  s, x\in {\bold R} \, .
$$
We also associate with $I_1$ the reflection $r_{I_1}$ given by
$$
r_{I_1}z = \bar z
$$
where $\bar z$ is the complex conjugate of $z$.  It follows from
the definition  that
$\Lambda _{I_1}$ restricts to an orientation preserving diffeomorphisms of
$I_1$, $r_{I_1}$ restricts to an orientation reversing diffeomorphism of
$I_1$ onto $I_1^\prime $. \par
Then, if $I$ is an interval and we choose $g\in {\bold G}$ such that
$I=gI_1$ we may set 
$$
\Lambda _I = g\Lambda _{I_1}g^{-1}\, ,\qquad
r_I = gr_{I_1}g^{-1}\, .
$$
%The elements $\Lambda _I(s)$, $s\in {\bold R}$ and $r_I$ are well defined. 
Let $r$ be an orientation reversing isometry of $S^1$ with
$r^2 = 1$ (e.g. $r_{I_1}$).  The action of $r$ on $PSL(2, {\bold R})$ by
conjugation lifts to an action $\sigma _r$ on ${\bold G}$, therefore we
may consider the semidirect product of
${\bold G}\times _{\sigma _r}{\bold Z}_2$.   Since
${\bold G}\times _{\sigma _r}{\bold Z}_2$ is a covering of the group
generated by $PSL(2, {\bold R})$ and $r$,
${\bold G}\times _{\sigma _r}{\bold Z}_2$ acts on $S^1$. We call
(anti-)unitary a representation $U$ of
${\bold G}\times _{\sigma _r}{\bold Z}_2$ by operators on ${\Cal H}$ such
that $U(g)$ is unitary, resp. antiunitary, when $g$ is orientation
preserving, resp. orientation reversing. \par
Now we are ready to define a  conformal
precosheaf. \par
An irreducible 
conformal precosheaf ${\Cal A}$ of von Neumann
algebras on the intervals of $S^1$
is a map
$$
I\rightarrow {\Cal A}(I)
$$
from ${\Cal I}$ to the von Neumann algebras on a separable Hilbert space 
${\Cal H}$ that verifies the following properties:
\vskip .1in
\noindent
{\bf A. Isotony}.  If $I_1$, $I_2$ are intervals and
$I_1 \subset I_2$, then
$$
{\Cal A}(I_1) \subset {\Cal A}(I_2)\, .
$$

\vskip .1in
\noindent
{\bf B. Conformal invariance}.  There is a nontrivial unitary
representation $U$ of
${\bold G}$  on
${\Cal H}$ such that
$$
U(g){\Cal A}(I)U(g)^* = {\Cal A}(gI)\, , \qquad
g\in {\bold G}, \quad I\in {\Cal I} \, .
$$                      
\vskip .1in
\noindent
{\bf C. Positivity of the energy}.  The generator of the rotation subgroup
$U(R(\vartheta ) )$ is positive.

\vskip .1in
\noindent
{\bf D.  Locality}.  If $I_0$, $I$ are disjoint intervals then
${\Cal A}(I_0)$ and $A(I)$ commute.

The lattice symbol $\vee $ will denote `the von Neumann algebra generated
by'.
\vskip .1in               
\noindent
{\bf E. Existence of the vacuum}.  There exists a unit vector
$\Omega $ (vacuum vector) which is $U({\bold G})$-invariant and cyclic for
$\vee _{I\in {\Cal I}}{\Cal A}(I)$.
\vskip .1in               
\noindent
{\bf F. Irreducibility}.  The only
$U({\bold G})$-invariant vectors are the scalar multiples of $\Omega$.

\vskip .1in
\noindent
The term irreducibility is due to the fact (cf. Prop. 1.2 of [GL]) that
under the assumption of {\bf F} $\vee_{I\in {\Cal I}} A(I) =B({\Cal H})$.
\par
We have the following (cf. Prop. 1.1 of [GL]):
\proclaim{2.1 Proposition} Let ${\Cal A}$ be an irreducible 
conformal precosheaf. The
following hold:
\roster
\item"{(a)}"  Reeh-Schlieder theorem: $\Omega $ is cyclic and separating
for each von Neumann algebra ${\Cal A}(I)$, $I\in {\Cal I}$.
\item"{(b)}"  Bisognano-Wichmann property: $U$ extends to an
(anti-)unitary representation of ${\bold G}\times _{\sigma _r}{\bold Z}_2$
such that, for any $I\in {\Cal I}$,
$$
\align
U(\Lambda _I(2\pi t)) &= \Delta _I^{it} \, \\
U(r_I) &= J_I \,
\endalign
$$             
where $\Delta _I$, $J_I$ are the modular operator and the modular
conjugation associated with $({\Cal A}(I), \Omega )$.
For each $g\in {\bold G} \times _{\sigma _r} {\bold Z}_2$
$$
U(g){\Cal A}(I)U(g)^* = {\Cal A}(gI) \, .
$$
\endroster
\roster
\item"{(c)}"  Additivity: if a family of intervals $I_i$ covers the
interval $I$, then
$$
{\Cal A}(I) \subset \vee _i {\Cal A}(I_i)\, .
$$
%\item"{(d)}"  $U$ is indeed a representation of
%$PSL(2, {\bold R})$, i.e. $U(2\pi ) = 1$.
\item"{(d)}" Haag duality: ${\Cal A}(I)' = {\Cal A}(I')$.
\endroster
\endproclaim   

Assume  ${\Cal A}$ is an irreducible conformal precosheaf  as defined in above.
By Cor. B.2 of [GL], the only invariant vectors 
under the action of the one-parameter group $U(\Lambda _I(2\pi t))$
are the scalar multiples of $\Omega$. This fact and (b) of Prop. 2.1 are 
usually referred to as the action of modular group is {\it ergodic}
and {\it geometric} respectively.  It follows that
${\Cal A}(I)$ is a type $III_1$ factor for any $I\in {\Cal I}$
(cf. Prop. 1.2 of [GL]).\par 
A covariant {\it representation} $\pi $ of
${\Cal A}$ is a family of representations $\pi _I$ of the
von Neumann algebras ${\Cal A}(I)$, $I\in {\Cal I}$, on a
separable Hilbert space ${\Cal H}_\pi $ and a unitary representation
$U_\pi $ of the covering group ${\bold G}$ of $PSL(2, {\bold R})$ 
such that the following properties hold:
$$
\align                
I\subset \bar I \Rightarrow \pi _{\bar I} \mid _{{\Cal A}(I)}
= \pi _I \quad &\text{\rm (isotony)} \\
\text{\rm ad}U_\pi (g) \cdot \pi _I = \pi _{gI}\cdot
\text{\rm ad}U(g) &\text{\rm (covariance)}\, .
\endalign
$$
A covariant representation $\pi$ is called irreducible if 
 $\vee _{I\in {\Cal I}}\pi({\Cal A}(I)) = B({\Cal H}_\pi)$. By our definition
the irreducible conformal precosheaf
 is in fact an irreducible representation 
of itself 
and we will call this representation the {\it vacuum representation}.
Note that by Cor. B.2 of [GL], the vacuum representation is
the unique (up to unitary equivalence) 
irreducible covariant representation which contains an eigenvector
of the generator of the rotation group with lowest eigenvalue $0$.  
A unitary equivalence class of covariant 
representations of ${\Cal A}$ is called
a {\it superselection sector} (cf. Page 17 of [GL]). \par
The composition of two superselection sectors, similar to  
the Connes's fusion (cf . Chap. V of [W2]), 
can be defined (cf. \S IV.2 of [FG] ). The composition is manifestly
unitary, and this is one of the virtues of the
above formulation. \par

Now let $G$ the group as in the introduction
. Let $g=Lie(G),$ 
$g_{\Bbb C}:= g\otimes {\Bbb C}$.
Denote by $\hat g$ the 
affine Kac-Moody algebra (cf. P. 163 of [KW]) associated to $g_{\Bbb C}$. 
Recall $\hat g=g_{\Bbb C}\otimes {\Bbb C}[t,t^{-1}] \oplus {\Bbb C} c$,
where ${\Bbb C}c$ is the 1-dimensional center of $\hat g$.\par
Denote by $LG$ the group of smooth maps
$f: S^1 \mapsto G$ under pointwise multiplication. The
diffeomorphism group of the circle $\text{\rm Diff} S^1 $ is
naturally a subgroup of $\text{\rm Aut}(LG)$ with the action given by
reparametrization. In particular ${\bold G}$ acts on $LG$. 
We will be interested
in the projective unitary representations 
(cf. Chap. 9 of [PS]) $\pi$ of $LG$ that
are both irreducible and have positive energy. This implies that $\pi $
should extend to $LG\ltimes \text{\rm Rot}$ so that
the generator of the rotation group Rot is positive.
It follows from
Chap. 9 of [PS] or P. 490 of [W2] 
that for fixed level (a finite set of positive integers,
see the footnote on this page)
, there are only finite number of such
irreducible projective 
representations (cf. 4.3 for a list in the case $G=SU(N)$).        
These irreducible projective 
representations can be  obtained by exponentiating the 
integrable representations of $\hat g$ at the same level by Th. 4.4 of [GW].
We refer the reader to III.7 of [FG] for a summary about the 
properties of these representations. \par
Now let $H\subset G$ be a connected Lie subgroup. 
%Let $g=Lie(G), h=Lie(H)$ and
%$g_{\Bbb C}:= g\otimes {\Bbb C}, h_{\Bbb C}:=  h\otimes {\Bbb C}$.
Let $\pi^i$ be an irreducible
projective representations of $LG$ with positive energy at level $L$
\footnotemark\footnotetext{
When G is the direct product of simple groups,
$L$ is a multi-index, i.e., $L=(L_1,...,L_n)$, where $L_i\in \Bbb {N}$
corresponding to the level of the $i$-th simple group. To save some writing
we write the coset as $H\subset G_L$ or
as $H\subset G$ when the level is kept fixed in the question.} on Hilbert
space $H^i$. When  restricting to $LH$, $\pi^i$ is a projective representation
of $LH$ with positive energy. By the proposition on P. 484 of [W2], 
$\pi^i$ is a direct sum of irreducible  projective representations
of $LH$.
Suppose when restricting to $LH$, $H^i$ decomposes as:
$$
H^i = \sum_\alpha H_{i,\alpha} \otimes H_\alpha 
,$$ and  $\pi_\alpha$ are irreducible projective representations of $LH$ on
Hilbert space $H_\alpha$.  The set of $(i,\alpha)$ which appears in
the above decompositions will be denoted by $exp$. 
The above decomposition can also 
be obtained via exponentiation (cf. Th. 4.4 of 
[GW]) from a similar decomposition (cf. (1.6.2) of [KW])
of integrable representation of
$\hat g$ when restricting to $\hat h$ 
(the affine Kac-Moody algebra associated with $H$).

\par
We
shall use $\pi^0$ (resp. $\pi_0$) to denote the vacuum representation of
$LG$ (resp. $LH$) on $H^0$ (resp. $H_0$).  
This is the unique projective   representation of
$LG$ (resp. $LH$) which contains a nonzero vector, unique up to
multiplication by a nonzero scalar, with the property that it is an eigenvector
of the generator of rotation group with eigenvalue 0. Such vectors will
be called vacuum vectors.  By Th. 3.2 of [FG] (also cf. \S17 of [W2]) 
the vacuume representation $\pi^0$ of  $LG$  gives rise to to 
an irreducible conformal precosheaf by the map
$I\in {\Cal I}\rightarrow  \pi^0 (L_IG)''$ on $H^0$. Similarly
$I\in {\Cal I}\rightarrow  \pi_0 (L_IH)''$ on $H_0$ is  
an irreducible conformal precosheaf on $H_0$.\par

Let $\Omega$ (resp. $\Omega_0$)  be a vacuum vector
in $\pi^0$ (resp.$\pi_0$) and assume
$$
\Omega = \Omega_{0,0} \otimes \Omega_0
$$ with $\Omega_{0,0} \in H_{0,0}$.
We shall always assume that $H\subset G_L$ is not a conformal inclusion
since
the case of conformal inclusions has been considered in [X1]. For the
definition of conformal inclusion, we refer the reader to P. 210 of [KW].
\par  
%According to [KW], this is equivalent to the fact that 
%$H_{i,\alpha}$ in the above decompositions are infinite dimensional.
%\par
\proclaim{Lemma 2.1}
(1)  $\pi^0(L_IH)''$ is  strongly
additive , i.e., if $I_1, I_2$ are the connected components
of the interval $I$ with one internal point removed, then:
$$
 \pi^0(L_IH)''= \pi^0(L_{I_1}H)'' \vee \pi^0(L_{I_2}H)'';
$$  \par
(2)
$\pi^0(L_IH)' \cap \pi^0(L_{I'}H)' = \pi^0(LH)'$.
\endproclaim
\demo{Proof}
Ad (1): Let $P_0$ be the
projection from $H^0$ onto $\Omega_{0,0}\otimes H_{0}$. 
%It follows from
%Th. 3.2 of [FG] or Prop. 1.2 of [GL] that for any
%interval $J$, $\pi^0(L_JG)''$ is a hyperfinite
%type $III_1$ factor. 
Since the action of the modular group of $\pi^0(L_JG)''$ with respect
to $\Omega$ is geometric and ergodic, it fixes globally $ \pi^0(L_JH)''$,
and  by Takesaki's theorem (cf. [MT] or P. 495 of [W2]),
$ \pi^0(L_JH)''$ is a factor.   So the map
$$
x\in \pi^0(L_JH)''\rightarrow xP_0
$$
is a $*$-isomorphism (cf. P. 492 of [W2]). 
Hence to prove (1) we just need to show 
$$
 \pi_0(L_{I}H)'' = \pi_0(L_{I_1}H)'' \vee \pi_0(L_{I_2}H)''. 
$$
By Haag duality in Prop. 2.1  it is enough to show
$$
 \pi_0(L_{I'}H)'' = \pi_0(L_{I_1'}H)'' \cap \pi_0(L_{I_2'}H)''. 
$$

%the lemma is trivial. Suppose $H$ is not a trivial group. 
By Reeh-Schlider theorem in Prop. 2.1,
the closed space spanned by $\pi^0(L_{J}H)\Omega$
is $P_0H^0$ for any interval $J$, and
by Takesaki's
theorem (cf. [MT] or (c) of Theorem on P. 495 of [W2]),
$$
\pi^0(L_{J}H)''= \{{\Bbb C} P_0 \}' \cap \pi^0(L_{J}G)''
,$$ so to prove (1) we just have to show that
$$
\pi^0(L_{I'}G)'' = \pi^0(L_{I_1'}G)'' \cap \pi^0(L_{I_2'}G)''. 
$$  It is sufficient to show the above in the case
when $G$ is simple, and this follows from Th. E of [W2]. \par  
Ad (2): 
Let $I_1=I, I_2, I_3, I_4$ be four  consecutive disjoint intervals on $S^1$
such that the closure of $I\cup I_2\cup I_3\cup I_4$ 
is $S^1$. Let $J_1=I_1', J_3=I_3'$. For any $a\in LH$,
since $H$ is connected, we can always choose $a_1\in L_{J_1}H$ such that
$a_1|\tilde I_1 = e, a_1|\tilde I_3= a|\tilde I_3$, where $e$ is the
identity of $G$,  
$I_1\subset \tilde I_1, I_3\subset \tilde I_3,$ and 
the closure of  $\tilde I_1$ and the closure of  $\tilde I_3$ are 
disjoint. Let $a_2 =  a a_1^{-1}$, then $a_2 \in L_{J_3}H$, and
$a_2 a_1 = a$. Hence $L_{J_1}H$ and $L_{J_3}H$ generate $LH$ algebraically,
and 
$$
\pi^0(L_{J_1}H)'' \vee \pi^0(L_{J_3}H)'' =\pi^0(LH)''.
$$
By (1)
$$
\pi^0(L_{I_2}H)'' \vee \pi^0(L_{I_{3}}H)'' \vee\pi^0(L_{I_4}H)'' = 
\pi^0(L_{J_1}H)'';
$$
$$
\pi^0(L_{I_1}H)'' \vee \pi^0(L_{I_{2}}H)'' \vee\pi^0(L_{I_4}H)'' = 
\pi^0(L_{J_3}H)'';
$$
Hence 
$$ 
\pi^0(L_{J_1}H)'' \vee \pi^0(L_{J_3}H)'' \subset
\pi^0(L_IH)'' \vee \pi^0(L_{I'}H)'',
$$ and
$$
\pi^0(LH)''= \pi^0(L_IH)'' \vee \pi^0(L_{I'}H)''.
$$
By taking the commutant of the above equality, the proof
of (2) is complete.
\enddemo 
\qed
\par
For each interval $I\subset S^1$, we define:
$$
A(I):=  \pi^0(L_IH)' \cap \pi^0(L_IG)''P
,$$ where $P$ is the projection from $H^0$ to closed subspace ${\Cal H}$
spanned by \par
$\vee_{J\in  {\Cal I}}\pi^0(L_JH)' \cap \pi^0(L_JG) \Omega$
.   Note that 
$$
\pi^0(L_IH)' \cap \pi^0(L_IG)'' \supset  \pi^0(LH)' \cap \pi^0(L_IG)''
.$$ On the other hand if $x\in  \pi^0(L_IH)' \cap \pi^0(L_IG)''$, then
$x\in \pi^0(L_IH)' \cap \pi^0(L_{I'}H)'$, but 
$\pi^0(L_IH)' \cap \pi^0(L_{I'}H)' = \pi^0(LH)'$ by (2) of lemma 2.1, so
$x\in\pi^0(LH)'$. This shows that
$$
\pi^0(L_IH)' \cap \pi^0(L_IG)'' =  \pi^0(LH)' \cap \pi^0(L_IG)''
.$$  It follows that if $x\in \pi^0(L_IH)' \cap \pi^0(L_IG)''$, then
as an operator on $H^0$, it takes the form $\oplus_{0,\alpha} B(H_{0,\alpha})
\otimes id_\alpha$, and  $x\Omega \in H_{0,0}\otimes \Omega_0$, 
so we have ${\Cal H}\subset H_{0,0}\otimes \Omega_0$.
\proclaim{Proposition 2.2} The map $I\in {\Cal I} \rightarrow
A(I)$ as defined above is an irreducible conformal precosheaf.
\endproclaim
\demo{Proof}
We have to check conditions {\bf A} to {\bf F}. {\bf A} follows from 
$$
\pi^0(L_IH)' \cap \pi^0(L_IG)'' =  \pi^0(LH)' \cap \pi^0(L_IG)''
$$ which
is proved above,
 {\bf B} and {\bf C} 
follows from a lemma  on P. 485 of [W2] except that we need to
show that the action of  ${\bold G}$ is nontrivial. 
Denote by $U(t)$ the action of the rotation group on $H_{0,0}$
and assume  $U(t)=\exp (2\pi iL_0^{g,h}t)$, where
$ L_0^{g,h}$ is 
the positive self-adjoint operator. Fix $\tau\in {\Bbb C}$ with $Im \tau>0$.
Then 
$$b^0_0(\tau)
:=\exp((1/24) 2 \pi i \tau (z_m -\dot z_{\dot m})) tr_{H_{0,0}} 
\exp(2\pi i \tau L_0^{g,h})
$$ is a branching function (cf. 3.2.7 of [KW]),
where $ (z_m -\dot z_{\dot m})$ is a nonnegative number defined as 1.4.2 of
[KW].  We will not need the detailed expression of
the branching function in our argument,
all we need to know is 
if  $H\subset G$ is not conformal, then $z_m -\dot z_{\dot m}>0$
(cf. P. 210 of [KW]). Let us show that if 
$H\subset G$ is not conformal which is assumed throughout this paper,
then  $ L_0^{g,h}\neq 0.$ 
If $ L_0^{g,h}=0,$ then 
$b^0_0(\tau)= \exp((1/24) 2 \pi i \tau (z_m -\dot z_{\dot m})) dim(H_{0,0})$
.  Since $b^0_0(\tau)$ is well defined for $Im\tau>0$ (cf. P. 170 of [KW]),
we must have $dim(H_{0,0})<\infty,$
but by 
\footnotemark\footnotetext{ The branching function here corresponds to
$b^0_0$ in the notations of [KW] on P. 187, where $0$ always denotes the
vacuum representations. The assumption of Th. B of [KW] follows from the
definition 2.5.4 of [KW] when $\Lambda =0, \lambda=0$.} (a) of Th. B of [KW],
$b^0_0(\tau) \sim b(0,0) \exp (\frac{\pi i (z_m -\dot z_{\dot m})}{12\tau})
$ as $\tau\rightarrow 0$ where $b(0,0)>0$. 
  This is a contradiction, and  shows that
 $ L_0^{g,h}\neq 0.$ 
So the action of ${\bold G}$ is nontrivial. 
{\bf F} follows from the uniqueness of the 
vacuum $\Omega$ (up to multiplication by a non-zero scalar) 
for $LG$. {\bf D} and {\bf E} follow from the definitions.
\enddemo \hfill \qed
\par
The irreducible conformal precosheaf
as in Prop. 2.2 is defined to be the irreducible conformal precosheaf
of the coset $H\subset G_L$. 
Note that when $H=\{e \}$ is a trivial subgroup
($e$ denotes the identity element in $G$), the irreducible 
conformal precosheaf
defined above coincides with the one defined in III.8 of [FG]. 
%It follows from
%Th. 3.2 of [FG] or Prop. 1.2 of [GL] that $\pi^0(L_IG)''$ is a hyperfinite
%type $III_1$ factor. 
Since the action of the modular group
of $\pi^0(L_IG)''$ with respect to $\Omega$ is geometric and ergodic, 
it fixes globally $ \pi^0(L_IH)''$,  hence
$ \pi^0(L_IG)'' \cap \pi^0(L_IH)'$, and by Takesaki's theorem (cf. [MT] or
P. 495 of [W2]),
$ \pi^0(L_IG)'' \cap \pi^0(L_IH)'$ is a   
factor. It follows that the map
$$ 
x\in  \pi^0(L_IG)'' \cap \pi^0(L_IH)' \rightarrow
xP  \in  \pi^0(L_IG)'' \cap \pi^0(L_IH)' P
$$ 
is a $*$ isomorphism (cf. P. 492 of [W2]), 
and can be implemented by a unitary $U_1: H^0 \rightarrow
{\Cal H}$, i.e., $ x = U_1^* xP U_1$, since 
$\pi^0(L_IG)'' \cap \pi^0(L_IH)' P$ is a type $III_1$ factor by Prop. 2.2
and the remarks after Prop. 2.1.
\par
Let us define a class of covariant representations of $A(I)$
coming from the decompositions of irreducible
projective representations $\pi^i$ of $LG$
with respect to $LH$. By the remarks after Th. B on P. 502 of [W2], for
any fixed interval $I$, 
there exists a unitary map $U: H^i\rightarrow H^0$ such that
$$
\pi^i(a) = U^* \pi^0(a) U, \forall a\in L_IG. 
$$ 
For $y=xP \in \pi^0(L_IH)' \cap \pi^0(L_IG)'P$, we define
$$
\pi^i (y) = U^*   U_1^* y U_1  U
.$$ This gives a factor representation
of $A(I)$. Let $P_{i,\alpha}$ be a projection from $H^i$ to
a subspace $H_{i,\alpha} \otimes \Omega_\alpha$ where $\Omega_\alpha$
is a unit vector in $H_\alpha$. Then
$$
y \in A(I) \rightarrow
\ \pi_{i,\alpha}(y) := \pi^i (y) P_{i,\alpha}
$$
is a subrepresentation of the factor representation
$\pi^i$, and so the map above is 
a $*$-isomorphism (cf. P. 492 of [W2]). 
Denote by $\pi^i(g)$ the action of $g\in {\bold G}$ on $H^i$. By the 
 lemma on P. 485 of [W2], $\pi^i(g)$ can be written as 
$\pi^i(g) = \oplus_\alpha \pi_{i,\alpha} (g) \otimes \pi_\alpha (g).$
One checks by using the definitions that the representations 
$ \pi_{i,\alpha}$ of $A(I)$ and   the representations 
$ \pi_{i,\alpha}$ of ${\bold G}$ satisfy the covariance condition,
and so  $ \pi_{i,\alpha}$ are 
covariant representations of $A(I)$. The study of these
representations is the main purpose of this paper.\par
By the same argument as in the proof of (1) of lemma 2.1  
one can show that $A(I)$ as
in Prop. 2.2 is strongly
additive , i.e., if $I_1, I_2$ are the connected components
of the interval $I$ with one internal point removed, then:
$$
A(I) = A(I_1) \vee A(I_2)
.$$  In fact, by Haag duality in Prop. 2.1, it is sufficient to show
that
$$
A(I') = A(I_1') \cap A(I_2')
,$$ which  is equivalent to
$$
\pi^0(L_{I'}H)'\cap \pi^0(L_{I'}G)'' =
(\pi^0(L_{I_1'}H)'
\cap \pi^0(L_{I_1'}G)'') \cap (\pi^0(L_{I_2'}H)'\cap \pi^0(L_{I_2'}G)'')
$$  by the paragraph after Prop. 2.2.
Let $P$ be the projection defined before Prop. 2.2.
By Reeh-Schlider theorem in Prop. 2.1
the closed space spanned by $\pi^0(L_{I'}H)'\cap \pi^0(L_{I'}G)''\Omega$
is $PH^0$ for any interval $I$, and
 by Takesaki's
theorem (cf. [MT] or (c) of Theorem on P. 495 of [W2]),
$$
\pi^0(L_{I'}H)'\cap \pi^0(L_{I'}G)'' = \{{\Bbb C} P \}' \cap \pi^0(L_{I'}G)''
,$$ so we just have to show that
$$
\pi^0(L_{I'}G)'' = \pi^0(L_{I_1'}G)'' \cap \pi^0(L_{I_2'}G)''.
$$  It is sufficient to show the above equation in the case
$G$ is simple, and in this case it follows from Th. E of [W2].\par
Note that 
the inclusion
$$
(\pi^0 ( L_{I}G)'' \cap\pi^0 (L_IH)') \vee \pi^0 (L_IH)'' \subset
\pi^0 ( L_{I}G)''
$$ is irreducible.
In fact 
$$
\align
& ((\pi^0 ( L_{I}G)'' \cap\pi^0 (L_IH)') \vee \pi^0 (L_IH)'')' \cap
\pi^0 ( L_{I}G)''\\
& =(\pi^0 ( L_{I}G)'' \cap\pi^0 (L_IH)')' \cap
 (\pi^0 ( L_{I}G)'' \cap\pi^0 (L_IH)') 
={\Bbb C},  
\endalign
$$ since $
\pi^0 ( L_{I}G)'' \cap\pi^0 (L_IH)'$ is a factor by the paragraph after 
Prop. 2.2. 
This fact can also be proved by using the fact that 
in $H_{0,0}$, the
vacuum representation of $A(I)$ appears once and only once. We shall
see in the next subsection that $\pi_{0,0}$ is in fact  the vacuum 
representation under general conditions.                 
\subheading {2.2  $\pi_{0,0}$ is the vacuum Representation}
As a first step in the study of representations $ \pi_{i,\alpha}$
we show that the representation of $ \pi_{0,0}$  on $H_{0,0}$ is 
the vacuum representation. This is equivalent to 
${\Cal H}= H_{0,0}\otimes \Omega_0$ by definition.
\proclaim{Theorem 2.3}
Suppose $H\subset G$,   and $H$ is   
simply connected.  Then
${\Cal H}= H_{0,0}\otimes \Omega_0$. Hence  $ \pi_{0,0}$ is the
vacuum representation.
\endproclaim
The idea of the proof is to use smeared vertex operators. From the
proof one can also see the close relation between our 
irreducible conformal precosheaf
 and
the definition of coset vertex operator algebra in \S5 of [FZ]. In fact, the 
coset vertex operator algebra in \S5  of [FZ] 
can be thought as ``germs'' of ours. \par
Let $g$ (resp. $h$) be the Lie
algebra of $G$ (resp. $H$).   
Choose a basis $e_\alpha,e_{-\alpha},h_\alpha$ in
$g_{\Bbb C}:=g\otimes {\Bbb C}$ with
$\alpha$ ranging over the set of roots as in \S 2.5 of [PS].
Let $X_\alpha :=e_\alpha + e_{-\alpha}, Y_\alpha:= i(e_\alpha -e_{-\alpha})$.
Denote by $\hat g$ the 
affine Kac-Moody algebra (cf. P. 163 of [KW]) associated to $g_{\Bbb C}$. 
Note $\hat g=g_{\Bbb C}\otimes {\Bbb C}[t,t^{-1}] \oplus {\Bbb C} c$,
where ${\Bbb C}c$ is the 1-dimensional center of $\hat g$. 
For $X\in g$,  Define $X(n):= X \otimes t^n$,  $X(z):= \sum_n X(n) z^{-n-1}$
as on Page 312 of [KT] and $
X^+(z):= \sum_{n<0} X(n) z^{-n-1},  X^-(z):= \sum_{n\geq 0} X(n) z^{-n-1}$. 
\par
Let $\pi^0$ be the vacuum representation of $LG$ on $H^0$ with vacuum
vector $\Omega$.  Let $D$ be the generator of the
action of the rotation group on $H^0$. 
For $\xi\in H^0$, we define $||x||_s = ||(1+D)^s x||, s\in {\Bbb R}$.
$H^{0,s}:=\{x\in H^0
|\ ||x||_s < \infty \}$ and $H^0_\infty = \cap_{s\in {\Bbb R}} H^{0,s}$.
Note that when $s\geq 0$, $H^{0,s}$ is a complete space under the norm 
$||.||_s$.
$H^0_0$ will denote the finite linear sum of the eigenvectors of $D$.
Clearly $H^0_0 \subset H^0_\infty$. The elements of $H^0_0$ (resp.
$H^0_\infty$) will be called {\it finite energy vectors} (resp. 
{\it smooth
vectors}). The eigenvalue of $D$ is sometimes referred to as 
energy. \par
Let us recall a few elementary facts about vertex operators which will
be used. See [FLM], [Kacv] or [D] for an introduction on vertex operator
algebras. Define $End(H^0_0)$ to be the space of all linear operators
(not necessarily bounded) from $H^0_0$ to $H^0_0$ and set
$$
End(H^0_0)[[z,z^{-1}]]:= \{ \sum_{n\in {\Bbb Z}} v_n z^n| v_n\in End(H^0_0)
\}.
$$
By the statement
on P. 154 of [FZ] which follows from Th. 2.4.1 of [FZ] there exists
a linear map
$$
\xi \in H^0_0 \rightarrow V(\xi,z)= \sum_{m\in {\Bbb Z}} \xi(m) z^{-m-1} \in
End(H^0_0)[[z,z^{-1}]] 
$$
with the following properties: \par
(1) $\xi(-1) \Omega = \xi$; \par
(2) If 
$$
\xi= 
X_{i_1}(-1)...
X_{i_t}(-1) \Omega
,$$ 
then 
$$
V(\xi,z)=  :X_{i_1}(z)...
X_{i_t}(z):
$$ where $:,:$ are normal ordered products (cf. (2.38), (2.39) of [D]).
\par
$V(\xi,z)$ is called a {\it vertex operator} of $\xi$. \par
Let $\psi \in H_{0,0} \otimes \Omega_0 $
be an eigenvector of $D$ with eigenvalue $n\in {\Bbb N}$. Then
$\psi$ takes the form:
$$
\psi = \sum_{} C_{i_1,...,i_t} X_{i_1}(-n_1)...
X_{i_t}(-n_t) \Omega
,$$ where the sum is finite , $n_i \geq 0$ and 
$C_{i_1,...,i_t} \in {\Bbb C}$. 
Since $g$ is semi-simple, $g=[g,g]$, so 
any $X(-n),n>1$ can be
expressed in terms of linear
combinations of the form 
$$
X(-1)Y(-1)...Z(-1).
$$  
Because $\Omega$ is the vacuum, 
$X\Omega =0, \forall X\in g$, so in  
$$ 
X_{i_1}(-n_1)...
X_{i_t}(-n_t) \Omega
,$$ if certain $X\in g$ appears,
we can always move $X$ to the right until it vanishes
when acting on $\Omega$.  For the above two 
reasons,  we can assume that
$n_1 =...=n_t = 1$. The vertex operator
$V(\psi, z)$ is then given by:
$$
V(\psi, z)=
\sum_{} C_{i_1,...,i_t} :X_{i_1}^{}(z)...
X_{i_t}(z):
,$$ where $:,:$ are normal ordered products by property (2) above.
\par 
Recall $ V(\psi, z)= \sum_m \psi(m) z^{-m-1}$. Define 
$$
V(m):= \psi(m+n-1)
$$
so we have
$ V(\psi, z)= \sum_m V(m)  z^{-m-n}$. 
This  expression for $ V(\psi, z)$ is in accordance with the 
convention of [KT].  Note that $V(-n)\Omega = 
\psi(-1)\Omega = \psi$ by property (1) above.
\par
Let $f=\sum_m f(m) z^m$ be a test
function with only a finite number of non-zero $f_m$.  Such $f$ will
be referred to as {\it finite energy functions}. Define
$$
||f||_s= \sum_{n\in {\Bbb Z}} (1+|m|)^s |f(m)|. 
$$
The 
{\it smeared vertex operator} $V(\psi,f)$ is defined to be:
$$
V(\psi,f) = \frac{1}{2\pi i}\int_{S^1} V( \psi, z) f dz = \sum_m f(m+n-1) V(m)
.$$ 
$V(\psi,f)$ is a well defined operator on $H^0_0$. 
Let $V(\psi,f)^{FA}$ be the formal adjoint of $V(\psi,f)$ on $H^0_0$.
It is defined by the equation
$$
\langle V(\psi,f)x, y\rangle = \langle x, V(\psi,f)^{FA}y \rangle
, \forall x,y\in H^0_0
$$
where $\langle, \rangle$ is the inner product on Hilbert space
$H^0$. 
When no confusion arises, we will write  $V(\psi,f)$ simply as
$V(f)$.  Similarly for $X\in g$, we define $X(f):=\sum_n X(n)f(n)$.\par
When the level is 1   
$H^0$ admits a fermionic representation (cf. \S13.3 of 
[PS] or I.6 of [W2]),
and we will denote the underlying Hilbert space by $F$. In fact, to
each simple component $G_i$ of $G$ there is a level 1 vacuum fermionic
representation of $LG_i$ on $F_i,i=1,...,m,$ and
$F=F_1\otimes F_2...\otimes F_m$.  
\proclaim{Lemma 1}
(1). 
Let $\xi \in H^0_0$, and $f$ is a finite energy function.
There exists 
 positive integer  $a$ and $c>0$ which are independent of  
$f$ and  $\xi$  such that 
$$
||V(\psi,f) \xi ||_s \leq c ||f||_{|s|+a} ||\xi||_{s+a}
;$$\par
(2). (1) is also true for $ V(\psi,f)^{FA}$ for the same constants
$c$ and $a$.
\endproclaim
\demo{Proof}
Note by definition 
$$
V(\psi,f) = \frac{1}{2\pi i}\int_{S^1} V( \psi, z) f dz = \sum_m f(m+n-1) V(m).
$$
If we can show (1) for the case 
$V(\psi,f) =V(l), \forall l$, then 
$$
\align
||V(\psi,f) \xi ||_s & = ||\sum_{l\in {\Bbb Z}} f(l+n-1) V(l) \xi||_s \\
&\leq c \sum_{l\in {\Bbb Z}} (1+|l|)^{|s|+a} |f(l+n-1)| \ ||\xi||_{s+a} \\
&\leq c (1+|n|)^{|s|+a} ||f||_{s+a} ||\xi||_{s+a}.
\endalign
$$
So we just need to show  (1) for the case 
$V(\psi,f) =V(l), \forall l$. Also note that $V(\psi,f)$ is linear
in $\psi$, so 
 it is sufficient to prove (1) in the case when 
$$
V(\psi, z) = :X_{1}^{}(z)...
X_{n}(z):
,$$ and $\xi$ is an eigenvector of $D$ with 
eigenvalue $\mu$. Similarly to prove (2) we just need to 
prove (2) in the case $V(\psi,f)^{FA} =V(l)^{FA}$, and 
$\psi, \xi$ are as above. Note that 
$$
:X_{1}(z)...
X_{n}(z):$$ is a summation of $2^n$ of expressions of the form
$ X_{i}^+...X_{j}^+X_{i'}^-...X_{j'}^-$. We will in fact prove
the ineqaulity in lemma 1 for such expressions. This will finish the proof of 
lemma 1 by definitions.  
We first  prove this for level 
$1$ and $G$ is simple, and the representation is on $F$. \par 

To avoid too many
subscripts we will denote the fermionic
creation or annilation operators simply by 
$a(m)$, since  we only need to use the fact that $a(m)$ increases the energy
by $m$, i.e., $[D, a(m)] = m a(m)$ on $H^0_0$ and its norm
is $1$ in our proof.     
Note in terms of $a(m)$,
$$
X(k)= \sum_{m>0} a(m-k) a(-m) - \sum_{m\geq 0} a(m) a(-m-k)
$$ when  acting on finite energy vectors,  cf. the expression
in a theorem on Page 486 of [W2].  Notice that
we have dropped all the subscripts here for simplicity. 
We will prove the lemma for $V(z) = X^+(z)...X^+(z)Y^-(z)...Y^-(z)$
where there are $n$ $X,Y$ and we have dropped the subscripts for
simplicity.
Then  $V(l)\xi$ is a sum of $2^n$ 
expressions of the form:
$$
\sum a(m_1)...a(m_{2n}) \xi
,$$ where the sum is over finite number of $m_1,...,m_n$'s, subject to certain
constraint. Let us now show by induction on $n$ that 
$$
0\leq |m_i| \leq e (1+|l| + \mu), i=1,...,2n, \tag 1
$$ where $e\geq 1$ depends only on $V(z)$.
When $n=1$, it is contained in a proof  on page 488  of [W2].
Assume that the statement is true for $k<n$. Suppose $V(z) = 
Y(z)X^-(z)$. Then:
 
$$
\sum_{k \geq 0} Y(l-k) X(k) = \sum_{k\geq 0, m>0} Y(l-k) a(m-k)a(-m)
- \sum_{k\geq 0, m\geq 0} Y(l-k) a(m) a(-m-k)
.$$ 
$ \sum_{k\geq 0, m>0} Y(l-k) a(m-k)a(-m) \xi$ is a sum of $2^{(n-1)}$ of
expressions of the form
$$
\sum_{k\geq 0,m>0, m_1,...,m_{2n-2}} a(m_1)...a(m_{2n-2}) a(m-k) a(-m) \xi
,$$ where $0\leq m_i \leq e' (1+|\mu -k| + |l-k|)$ by induction hypothesis,
and $e'\geq 1$ depends only on $Y(z)$.
Note that the above expression is nonzero only if $ 0< m \leq \mu,
0\leq k \leq \mu$. It follows that the expression is the sum of
$$
\sum_{m_1,...,m_{2n}} a(m_1)...a(m_{2n}) \xi
$$ with $0\leq |m_i| \leq e (1+|l| + \mu), i=1,...,2n$, where 
$e = 4 e'$. The same conclusion holds for 
$\sum_{k\geq 0,m>0, m_1,...,m_{2n-2}} a(m_1)...a(m_{2n-2}) a(m) a(-m-k) \xi$
.\par
Suppose $V(z) = 
X^+(z) Y(z)$. Then
$$\sum_{k < 0} X(k) Y(l-k)  = \sum_{k< 0, m>0}  a(m-k)a(-m) Y(l-k)
- \sum_{k< 0, m\geq 0}  a(m) a(-m-k) Y(l-k)
.$$ 
$ \sum_{k<  0, m>0}  a(m-k)a(-m) Y(l-k)\xi$ is a sum of $2^{(n-1)}$ of
expressions of the form
$$
\sum_{k<0,m>0, m_1,...,m_{2n-2}}  a(m-k) a(-m) a(m_1)...a(m_{2n-2})  \xi
,$$ where $0\leq m_i \leq e'' (1+|\mu | + |l-k|)$ by induction hypothesis,
and $e''\geq 1$ depends only on $Y(z)$.  
Note that the above expression is nonzero only if $ 0< m \leq \mu - l +k,
l-\mu < k < 0$.   It follows that the expression is the sum of
$$
\sum_{m_1,...,m_{2n}} a(m_1)...a(m_{2n}) \xi
$$ with $0\leq |m_i| \leq \tilde e (1+|l| + \mu), i=1,...,2n$, where 
$\tilde e = 4 e''$.
The same conclusion holds for 
$\sum_{k<0,m\geq 0, m_1,...,m_{2n-2}}  a(m) a(-m-k) a(m_1)...a(m_{2n-2})  \xi$
.  By induction (1) is proved.
Since the norm of $a(m)$ is $1$, it follows from (1) that:
$$
||V(l) \xi||_s \leq c (1+|\mu -l|)^s (1+\mu +|l|)^a ||\xi||
,$$ where $c,a$ are independent of $\xi$ and $l$.
Note $  (1+|\mu -l|)^s (1+\mu +|l|)^a \leq (1+\mu)^{s+a} (1+|l|)^{s+a} 
$ when $s\geq 0$.  When $s<0$, the function
$$
(\frac{1+x}{1+|x-l|})^{-s}
$$ has maximum $(1+l)^{-s}$ when $x\geq 0, l\geq 0$, and $1$ when
$x\geq 0, l< 0$. It follows that if $s<0$,
$$
(1+|\mu -l|)^s (1+\mu +|l|)^a \leq (1+\mu)^{s+a} (1+|l|)^{-s+a}
.$$
So we have:
$$
 ||V(l) \xi||_s \leq c (1+|l|)^{|s|+a} ||\xi||_{s+a}
.$$ \par 
Now let us consider the case when $G$ is simple and 
level $k>1$. On $F^{\otimes k}$, $  X_{i}^+...X_{j}^+X_{i'}^-...X_{j'}^-$
is a summation of $k^n$ of expressions of the form:
$A:=Y_1 \otimes Y_2 ...\otimes Y_l$, and each $Y_i$ is of the form
$  X_{p}^+...X_{q}^+X_{p'}^-...X_{q'}^-$.
Let $\xi_1 \otimes ...\otimes \xi_k$ be a vector in $F^{\otimes k}$
with $D\xi_i = \mu_i \xi_i,i=1,...,k$. Then:
$$
\align
& ||A(l)\xi_1 \otimes ...\otimes \xi_k||_s = ||\sum_{m_1,...,m_k, \sum_im_i =l}
Y_1(m_1)\xi_1 \otimes ...\otimes Y_k(m_k)\xi_k||_s \\
& \leq (1+|\mu_1+...+\mu_k -l|)^s \sum_{m_1 ,...,m_k, \sum_im_i =l} 
||Y_1(m_1)\xi_1 \otimes ...\otimes Y_k(m_k)\xi_k|| \\
& =   (1+|\mu_1+...+\mu_k -l|)^s \sum_{m_1 ,...,m_k, \sum_im_i =l, m_i\leq 
\mu_i} 
||Y_1(m_1)\xi_1 \otimes ...\otimes Y_k(m_k)\xi_k|| \\
& \leq (1+\mu_1+...+\mu_k)^s (1+|l|)^{|s|} \times  \\
&\sum_{m_1 ,...,m_k, 
\sum_im_i =l, m_i\leq \mu_i} c'
(1+|m_1|)^{b'}...(1+|m_k|)^{b'} ||\xi_1||_{a'}...
||\xi_k||_{a'} \\
& \leq c' ( 1+|l|)^{|s|+b'k+k}||\xi_1 \otimes ...\otimes \xi_k||_{s+b'k+k+a'k}
\endalign
$$ where in the second $\leq$ we use the result for $k=1$, and 
in the last $\leq$ we used $l- \sum_{j\neq i} \mu_j \leq m_i \leq \mu_i,
i=1,...,k.$  
The constants $a',b',c'$ above are independent of $l, \xi_1,...,\xi_k$.\par
Assume $G=G_1\times...\times G_m$ where $G_i,i=1,2,...,m$ are simple
factors of $G$. Let $F_i$ be the level 1 vacuum fermionic representation
of $LG_i, i=1,...,m.$  On $F_1^{k_1} \otimes ... \otimes F_m^{k_m}$,
 $  X_{i}^+...X_{j}^+X_{i'}^-...X_{j'}^-$
is still a finite sum  of expressions of the form:
$A:=Y_1 \otimes Y_2 ...\otimes Y_l$, and each $Y_i$ is of the form
$  X_{p}^+...X_{q}^+X_{p'}^-...X_{q'}^-$. Similar argument as in the
previous paragraph shows that the lemma is true for such expressions. 
Since   level $(k_1,...,k_m)$ vacuum representation of $L(G_1\times
G_2\times...\times G_m)$ 
appears as a direct summand of $F_1^{k_1} \otimes ... \otimes F_m^{k_m}$, 
we have proved the lemma. 
%(2) of the lemma is proved in similar way since by 
%definition $V(\psi,f)^{FA}$
%are expressed in terms of expressions similar to $V(\psi,f)$.
\enddemo 
\qed
\par
By (1) of lemma 1 the domain of $\pi^0(V(f))$ can
be extended to $H_\infty^0$, and  $\pi^0(V(f))$ can be defined for smooth $f$. 
%We claim that
%$H^0_0$ the adjoint $ \pi^0(V(f))^*$
When no confusion arises, we will denote 
$\pi^0(V(f))$ simply by $V(f).$
% and since 
%by (2) of lemma 1 the formal adjoint of $\pi^0(V(f))$ has a dense
%domain, it follows that the operator $\pi^0(V(f))$ is closable 
%and we will denote the closure of these operators with
%domain $H^0_\infty$ by the same symbol.  

%Also
%note that these operators are well defined on $H_\infty^0$
%for smooth $f$ by Lemma 1. 
Recall that an operator (not necessarily bounded) $a$ is called 
{\it affiliated with a von Neumann algebra M} if $U^* a U=a$ for any 
unitary $U\in M'$ (cf. P. 16 of [Dix]). 
Our goal is to show that  $V(f)$ is  
affiliated with certain von Neumann algebra (cf. Prop. 2.3 below).
%When no confusion arises, we will denote operators
%$\pi^0(V(f))$ simply by $V(f).$
\proclaim{Lemma 2}
Assume  $f$ is
a smooth function. Then: \par
%(1) The space of finite energy vectors $H^0_0$ is a core for $V(f)
%$; \par
%(2) Suppose $M$ is a von Neumann algebra  on $H^0$ with a generating
%self adjoint subset $S\subset M$,i.e., $S^*=S$ and the $C^*$ algebra
%generated by $S$ is dense in $M$ in strong topology. If  
%$$
%V(f) s x = s V(f)x, \forall s\in S, x\in H^0_0
%,$$ then $V(f)$ is affiliated with $M'$, $j=1,2$.\par
(1) Denote by $V(f)^*$ the adjoint of $V(f)$, then 
$H^0_\infty$ is in the domain of  $V(f)^*$. \par
(2)   The operator $V(f)$ is closable 
and (we will denote the closure of  $V(f)$ 
by the same symbol)
the space of finite energy vectors $H^0_0$ is a core for $V(f)
$; \par
(3) Suppose $M$ is a von Neumann algebra  on $H^0$ with a generating
self adjoint subset $S\subset M$,i.e., $S^*=S$ and the $C^*$ algebra
generated by $S$ is dense in $M$ in strong topology. If  
$$
V(f) s x = s V(f)x, \forall s\in S, x\in H^0_0
,$$ then $V(f)$ is affiliated with $M'$, $j=1,2$.\par
%(2) The space of finite energy vectors $H^0_0$ is a core for $V(f)
%$;
\endproclaim
\demo{Proof}
Ad (1): It is enough to show
that for any $y\in H^0_\infty$, there exists $w\in H^0$ such that
$$
\langle V(f)x, y\rangle = \langle x, w\rangle, \forall x\in H^0_\infty.
$$  
%First we claim it is enough to show the above equation for any
%$x\in H^0_0$. Suppose the above  equation is true for
%any $x\in H^0_0$. 
Note that if $x\in H^0_\infty$, then there exists a sequence of
finite energy vectors $x_n$ such that
$$
lim_{n\rightarrow\infty}||x_n-x||_a=0  
$$
where $a>0$ is as in lemma 1. By lemma 1
$$
lim_{n\rightarrow\infty}||V(f)x_n- V(f)x||_a=0.  
$$ 
So if 
$$
\langle V(f)x_n, y\rangle = \langle x_n, w\rangle, \forall x_n,
$$
then
$$
\langle V(f)x, y\rangle = \langle x, w\rangle, \forall x\in H^0_\infty.
$$
So we just have to prove the above equation for 
$x\in  H^0_0$. \par
First we consider the case that $f$ is a finite energy function.
 Since $y\in H^0_\infty$,
we can choose a  sequence $y_{m}$ such that
$$
lim_{m\rightarrow\infty}||y_m-y||_a=0  
$$  where $a>0$ is as in lemma 1.
By definition we have
$$
\langle V(f)x, y_m\rangle = \langle x,  V(f)^{FA}y_m\rangle,
$$ and by (2) of lemma 1 
$$
||V(f)^{FA}(y_m-y_{m'})||\leq c ||f||_a ||y_m-y_{m'}||_a.
$$
So
$\{V(f)^{FA}y_m\}_{m\geq 0}$ is a Cauchy sequence with a limit
defined to be $V(f)^{FA}y$.  We  can choose $w=V(f)^{FA}y$.
Also note that
$$
||V(f)^{FA}y||\leq c ||f||_a ||y||_a.
$$ 
\par
Now let  $f$ be a smooth function, and  choose a sequence 
$f_n$ of finite energy functions such that
$$
lim_{n\rightarrow\infty}||f_n-f||_a=0 
$$
where constant $a$ is as in lemma 1. 
%Similarly since $y\in H^0_\infty$,
%we can choose a double sequence $y_{m}$ such that
%$$
%lim_{m\rightarrow\infty}||y_m-y||_a=0.  
%$$  
By the proof in the finite energy function case  we have
$$
\langle V(f_n)x, y\rangle = \langle x,  V(f_n)^{FA}y\rangle.
$$  
By  lemma 1 and the note above
$$
||V(f_n-f)x ||\leq c ||f_n-f||_a ||x||_a, \ 
||V(f_n-f_{n'})^{FA}y||\leq c ||f_n-f_{n'}||_a ||y||_a.
$$
It follows that the sequnce $ V(f_n)^{FA}y,\forall n> 0$ is
a Cauchy sequence with a limit denoted by $w$ and 
$$
\langle V(f)x, y\rangle = \langle x, w\rangle, \forall x\in H^0_0.
$$  
\par
Ad (2): By (1) $V(f)$ is closable.  
Let $x\in H^0_\infty$. Then one can find $x_n\in H^0$ such that
$||x_n-x||_a \rightarrow 0$ with
$a>0$ as in  lemma 1, and by lemma 1 $V(f)x_n \rightarrow V(f)x$.
This shows $H^0_0$ is a core for $V(f)$. \par
Ad (3): 
Suppose
$y_n\rightarrow y\in M$ in the strong topology and 
$$
V(f) y_n x = y_n V(f)x, \forall s\in S, 
$$ and for all $x\in$ the domain of $V(f)$, it follows immediately that
$yx$ is in the domain of $V(f)$ and
$$
V(f) y x = y V(f)x, \forall s\in S, 
$$ and for all $x\in$ the domain of $V(f)$. \par 
Since 
$$
V(f) s x = s V(f)x, \forall s\in S, x\in H^0_0
,$$ and   $H^0_0$ is a core for $V_j(f)$, it follows that
$$
V(f) s x = s V(f)x, \forall s\in S, 
$$ and for all $x\in$ the domain of $V(f)$.
So 
$$
V(f) s_1s_2...s_n x = s_1s_2...s_n V(f)x, \forall s_i\in S, i=1,...,n 
$$ and for all $x\in$ the domain of $V(f)$ and finite $n$. 
(3) now follows from the definition.\par
%Note that the set $\{||f_n||_a, \forall n\geq 0\}$ is bounded, so
%$\{V(f_n)^{FA}(y_m)\}_{m\geq 0}$ is a Cauchy sequence with a limit
%defined to be $V(f_n)^{FA}y$. Note that we only use the fact that 
%$f_n$ is a finite energy function 
%$$
%||V(f_n)^{FA}y||\leq c ||f_n||_a ||y||_a,
%$$  
\enddemo
\hfill \qed
\par
Suppose $p^*=-p$ is a smooth test function. 
Assume that $X\in g$.
It follows from \S3 of [GW]
(also cf. P. 489 of [W2])  that
$X(p)$ is essentially skew-self adjoint with core $H^0_0$,
and    $X(p)$ maps  $H^0_\infty$ to  $H^0_\infty$.
\proclaim{Lemma 3}
Let
$X\in g$. Then: \par
(1) 
$$
[X(p), V(\psi,f)]x = 
\sum_{0\leq j \leq n}  V(X(j)\psi,  \frac{1}{j!}\frac{d^j p}{dz^j}  f (z)  )x
$$ for any smooth functions $p,f$ and $x\in H^0_\infty$, where
$n$ is the energy of $\psi$;\par
(2) 
If $\exp(tX(p)) H^0_0 \subset H^0_\infty, -1\leq t\leq 1$, and
$Sup\{||\exp(tX(p)) x||_s, -1\leq t\leq 1 \} < \infty$ for
any $x\in H^0_0, s>0$,  then
$$
\langle [\exp(X(p)), V(\psi,f)] x, y\rangle
 = \int_0^1 \langle \exp(tX(p))[X(p),V(\psi,f)] 
\exp((1-t)X(p)) x, y\rangle dt
,$$ for any smooth functions $p=-p^*,f$ and $x\in H^0_0, y\in H^0_0$.
\endproclaim
\demo{Proof}
Note for $x,y\in H^0_0$, 
$<X(k) V(\psi,z)x - V(\psi,z) X(k)x ,y>$ is a  polynomial in
$z,z^{-1}$.

We have (cf. for an example P. 327 of [KT] ): for $z\neq 0$,
$$
\align
<X(k) V(\psi,z)x - & V(\psi,z) X(k)x ,y>  = \frac{1}{2\pi i}\int_{C_z} dw w^k <
X(w) V(\psi,z) x, y> \\
&=  \frac{1}{2\pi i}\int_{C_z} dw w^k \sum_{m} (w-z)^{-m-1}<V(X(m)\psi,z)x,y>
\\
&=  \frac{1}{2\pi i}\int_{C_z} dw w^k \sum_{m\leq n} (w-z)^{-m-1}<
V(X(m)\psi,z)x,y>
\\
&= \sum_{0\leq j \leq n}
 \frac{1}{j!} \frac{d^j z^k}{dz^j}   <V(X(j) \psi, z)x,y> 
\endalign
$$ where $C_z$ is the boundary of a disk  centered at $z$ with radius $\frac
{1}{2}|z|$, and in the last equation we used the fact that
$$
\frac{1}{2\pi i}\int_{C_z} dw w^k (w-z)^{-j-1}  
=\frac{1}{j!} \frac{d^j z^k}{dz^j}   
$$
for $0\leq j \leq n$.
Since $H^0_0$ is dense,  we have:
$$
\align
[X(k), V(\psi,f)]x & = \sum_{0\leq j \leq n}  \frac{1}{2\pi i} \int_{S^1}
 \frac{1}{j!} \frac{d^j z^k}{dz^j}  f (z) V(X(j)\psi, z) \\
& =  \sum_{0\leq j \leq n}  V(X(j)\psi,  \frac{1}{j!} \frac{d^j z^k}{dz^j}  
f (z)  )x
\endalign
$$ which is  true for any finite energy function $f$ 
 and
$x\in H^0_0$, and so it is true for any smooth function
$f$ and $x\in H^0_\infty $ by using approximation and lemma 1 .
%By the same argument using lemma 1, we have:
Let $p=\sum_k p(k)z^k$ be a finite energy
function and $x\in H^0_\infty$. By definition
$$
X(p)= \sum_k p(k) X(k),
$$ and so (remember $p=\sum_k p(k)z^k$)
$$
[X(p), V(\psi,f)]x =
\sum_{0\leq j \leq n}  V(X(j)\psi,  \frac{1}{j!} 
\frac{d^j p}{dz^j}  f (z)  )x.
$$ 
By lemma 1 the above is also true for any smooth function $p$ since
we can always choose a sequence of functions $p_m$,  each $p_m$ is 
a finite energy function, and $||p_m-p||_s\rightarrow 0$ as
$m\rightarrow \infty$
for $s$ greater than a given number which may depend on 
$X, \psi, X(j)\psi, j=0,...,n$.\par 
Ad (2): Note $H^0_\infty$ is a subset of $C^\infty$ vectors of $X(p)$ and 
$V(f):=V(\psi,f)$.   
Let us  check that 
the map 
$$
s\in [0,1]\rightarrow A(s):=\langle \exp(sX(p))V(f)\exp((1-s)X(p))x, 
y\rangle
$$ is a differentiable
function  with continuous derivative 
$$
B(s):=\langle \exp(sX(p))[X(p),V(f)] \exp((1-s)X(p)) x, y\rangle. 
$$
Define 
$$
C(s,t):= \langle \exp(sX(p))V(f)\exp((1-t)X(p))x, y\rangle
, (s,t)\in [0,1]\times [0,1].$$ 
We shall repeatedly use the following elementary fact 
about $C^\infty$ vectors (cf. P. 488 of [W2]): if
$\xi$ is a $C^\infty$ vector of $X(p)$, i.e., 
$\xi$ is in the domain of $X(p)^n, \forall n\geq 1$, then
the function
$$
u\in {\Bbb R}\rightarrow \langle \exp(uX(p)) \xi, \eta \rangle
$$
is a smooth function of $u$ for any $\eta\in H^0$. \par
Since $V(f)\exp((1-t)X(p))x \in H^0_\infty$ is a subset
of  $C^\infty$ vectors of $X(p)$, it follows that
$C(s,t)$ is a smooth function of $s$ for fixed $t$. Also note
$$
\align
C(s,t)&= \langle \exp(sX(p))V(f)\exp((1-t)X(p))x, y\rangle \\
&= \langle \exp((1-t)X(p))x, V(f)^*\exp(-sX(p))y\rangle
\endalign
$$ where we have used (1) of lemma 2 and our assumption that
$\exp(-sX(p))y  \in H^0_\infty, 0\leq s \leq 1$.  So
$C(s,t)$ is a smooth function of $t$ for fixed $s$. We have  
( use (1) of lemma 2 when computing derivatives with respect to $t$)  
the following partial derivatives:
$$
\align
C_s(s,t) &= \langle \exp(sX(p)) X(p)V(f)\exp((1-t)X(p))x, y\rangle \\
C_{ss}(s,t) &= \langle \exp(sX(p)) X(p)^2V(f)\exp((1-t)X(p))x, y\rangle \\ 
C_t(s,t) &= \langle \exp(sX(p)) V(f) (-X(p))\exp((1-t)X(p))x, y\rangle \\
C_{tt}(s,t) &= \langle \exp(sX(p)) V(f) X(p)^2\exp((1-t)X(p))x, y\rangle \\
C_{st}(s,t) = C_{ts}(s,t) &= 
\langle \exp(sX(p)) X(p)V(f) (-X(p))\exp((1-t)X(p))x, y\rangle 
\endalign
$$
Note that all the derivatives above are smooth functions of
one variable when the other variable is fixed. 
We have
$$
\align
|C_{ss}(s,t)|& \leq ||\exp(sX(p)) X(p)^2V(f)\exp((1-t)X(p))x|| \ ||y||\\
& \leq  ||X(p)^2V(f)\exp((1-t)X(p))x|| \ ||y||\\
&\leq C_1||p||^2_{a_1} ||f||_{a_2}  ||\exp((1-t)X(p))x||_{a_3}||y||\\
&\leq C_2
\endalign
$$
where we used lemma 1 in the third $\leq$, 
the assumption in the last $\leq$,
and 
$$
C_1, C_2, a_1, a_2, a_3
$$
are independent of $s,t$ by lemma 1 and the assumption. We can
obtain similar estimates for other second partial derivatives, and so
there exists a constant $C'>0$ such that
$$
Sup\{|C_{ss}(s,t)|, |C_{st}(s,t)= C_{st}(t,s)|, |C_{tt}(s,t)|, 
\forall (s,t)\in [0,1]\times [0,1] \} \leq C'.
$$
By using the uniform bound for second partial derivatives and Taylor's theorem
in calculus, we have
$$
\align
&A(s+\Delta s) -A(s)
= C(s+ \Delta s, s+ \Delta s) -C(s,s) \\
&= C(s+ \Delta s, s+ \Delta s) -C(s+ \Delta s,s)+ 
C(s+ \Delta s, s) -C(s,s) \\
&= 
C_t(s+\Delta s, s) \Delta s +  \frac{1}{2}C_{tt} 
(s+\Delta s, \theta_1) (\Delta s)^2 + C_s (s,s)\Delta s+ 
\frac{1}{2} C_{ss} (\theta_2,s)(\Delta s)^2 \\
&= (C_t(s, s) + \frac{1}{2} C_{ts}(\theta_3, s) \Delta s) \Delta s
+ \frac{1}{2}C_{tt} 
(s+\Delta s, \theta_1) (\Delta s)^2 + \\
&C_s (s,s)\Delta s+ 
\frac{1}{2} C_{ss} (\theta_2,s)(\Delta s)^2 \\
&= B(s) \Delta s + O((\Delta s)^2)
\endalign
$$  
where  $\theta_i, i=1,2,3$ are
between $s$ and $s+\Delta s$, 
$$|O((\Delta s)^2)|\leq \frac{3C'}{2}|(\Delta s)^2|,
$$ and 
we have used 
$$
B(s) = C_t(s, s) + C_s (s,s)
$$
which follows from definitions. It follows immediately that
the derivative of $A(s)$ is $B(s)$ on [0,1]. A similar 
elementary exercise in calculus as above shows that
$B(s)$ is continuous on [0,1]. (2) now follows by the Fundamental Theorem
of Calculus. 
\enddemo
\hfill \qed
\par
Let $T$ be the maximal torus of $G$ and $\top =Lie(T)$. 
By \S13.3 of [PS], the level 1 vacuum
representation of $LT$ on Hilbert space $F$, 
is also an irreducible representation of $LG$.
Denote by $\pi$ the representation of $LG$ on $F$.
\proclaim{Lemma 4}
(1) Let $k\in {\Bbb N}$ and $x \in F_0,  $ where $F_0$ denotes
the set of finite energy vectors and
$v=\exp(w(p))$ with $w\in \top, p=-p^*,||p||_{k+1}< M $.
 Then there exists a constant $C$
which only depends on $w, M, k$ and $x \in F_0  $ such that
$$
||\pi (v)x||_k \leq C
;$$ \par
(2) Let $u=\exp(X(p))$ where $X=X_\alpha$ or
$X= Y_\alpha$ and $p=-p^*, ||p||_{k+1} <M$. 
Then  then there exists a constant $C'$
which only depends on $X, M, k$ and $x \in F_0  $ such that
$$
||\pi(u)x||_k \leq C'
;$$ \par
(3) Let $u=\exp(X(p))$ where $X=X_\alpha$ or
$X= Y_\alpha$ and $p=-p^*, ||p||_{k+1} <M$. Denote by $\pi^0$
the vacuum representation of $LG$ on $H^0$. 
Then  then there exists a constant $C''$
which only depends on $X, M, k$ and $x \in H^0_0  $ such that
$$
||\pi^0(u)x||_k \leq C'' 
.$$
\endproclaim
\demo{Proof}
Ad (1): The basic idea is contained in Prop. 9.5.15 of [PS] and
we will recall the  notations and facts in 9.5 [PS]. We can 
write $LT \simeq Hom(S^1,T) \times T\times
V$, where $T$ is the subgroup of constant loops, and $V$ is the vector
space of maps $a :S^1 \rightarrow \top$ with integral $0$, which is 
regarded as a subgroup of $LT$ by the exponential map
$a\rightarrow \exp(ia)$.
The identity component of the central extension of
$LT$ in our case is canonically a product $T\times \tilde V$, where $\tilde V$
is the Heisenberg group associated to a skew form $S$ defined  on P. 
63 of [PS]. 
Write $V\otimes {\Bbb C}= A\oplus \bar A$,  where $A$ is spanned
by $z^k \top \otimes {\Bbb C}$ for $k>0$.  
For $a\in V\otimes {\Bbb C}$, let $a=\sum_n a_n z^n, a_n\in 
\top \otimes {\Bbb C}$ be 
its Fourier series. We define 
$$
||a||_s := \sum_n (1+|n|)^s |a_n|, s\in {\Bbb R},
$$  
and $|.|$ is the norm on $\top\otimes {\Bbb C}$ 
induced from the restriction of 
Killing form on $\top$. The Hermitian form
$\langle, \rangle$ on $A$ defined by 
$$
\langle  a, a' \rangle
= -2i S (\bar \xi, \eta) 
$$
is positive definite, where the  skew form $S$ is  defined on P. 63 of 
[PS].  The only property we need about $S$ is 
$$
|S(a, a')| \leq ||a ||_0 ||a'||_1, \forall a,a'\in V
$$ which follows from its definition. \par
The Hilbert space $F$
is the completion of the symmetric algebra $S(A)$ with 
respect to Hermitian form above, which is extended from $A$ to $S(A)$ by the
formula
$$
\langle  a_1a_2...a_n, a'_1a'_2...a'_n \rangle
= \sum \langle a_1, a'_{i_1}\rangle...\langle a_n, a'_{i_n}\rangle
$$
where the sum is over all permutations $\{ i_1,...,i_n\}$ of
$\{ 1,...,n\}$.  Note that for any $a\in A$, 
$$
e^a:=\sum_{n\geq 0} \frac{a^n}{n!}
$$
belongs to $F$. The vacuum vector in $F$ is denoted by $1$. \par 
As in the proof of  Prop. 9.5.15 of [PS], it is sufficient to
prove (1) for the case when $v\in \tilde{V}$ and $x$ is the
vacuum vector. Note that $v=\exp(w(p))$ is identified with
$ \tilde w(p)= i pw \in \tilde{V}$ under
the isomorphism 
$LT \simeq Hom(S^1,T) \times T\times
V$ above. 

The action
of $v$ on vacuum vector $1$ is given by (cf. P. 194 of [PS]):
$$
v.1 = e^{-\frac{1}{2} \langle a,a\rangle} e^a
,$$ where $\tilde w(p)=a+ \bar a$, and  $a(z) = \sum_{i> 0} a_i z^{i}$. 
Let 
$$ 
a^{(s)}(z):= \sum_{i>0} i^s a_i z^{i},
 s\in {\Bbb N}
.$$
Note that $||a^{(s)}||_0 \leq ||a||_s \leq ||p||_s |w|$.  
We have:
$$
D^k a^n= \sum_{s_1\geq 0,... s_n\geq 0,s_1+...+s_n =k}
\frac{k!}{s_1!...s_n!}a^{(s_1)}... a^{(s_n)}  
,$$ and so
$$
||D^k a^n||\leq \sum_{s_1\geq 0,... s_n\geq 0,s_1+...+s_n =k}
\frac{k!}{s_1!...s_n!}||a^{(s_1)}... a^{(s_n)}||  
.$$
Note that for $ 0\leq s_1,t_1 \leq k$, 
$$
\align
|\langle a^{(s_1)}, a^{(t_1)} \rangle|
& = |2S(a^{(s_1)}, \overline{a^{(t_1)}})| 
\leq 2 || a^{(s_1)}||_0 ||a^{(t_1)}||_1 \\
&\leq  2 || a||_{s_1} ||a||_{t_1+1} \leq 2 ||p||^2_{k+1} |w|^2,
\endalign
$$ hence 
$$
\align
||a^{(s_1)}... a^{(s_n)}||^2 & = 
\langle a^{(s_1)}... a^{(s_n)},  a^{(s_1)}
... a^{(s_n)} \rangle \\
& \leq n! 2^n |w|^{2n}  ||p||^{2n}_{k+1}
\endalign  
$$
where we used the definition of $\langle, \rangle$ on 
$F$ as completion of $S(A)$ . 
So  
$$
||D^k a^n|| \leq  (n!)^{\frac{1}{2}} ||p||_{k+1}^{n} n^{k} (\sqrt{2}|w|)^{n}
,$$ and
$$
\align
||D^k v.1|| & \leq  e^{-\frac{1}{2} \langle a,a\rangle} \sum_{n=1}^\infty
\frac{1}{(n!)^{\frac{1}{2}}} ||p||_{k+1}^{n} n^{k} (\sqrt{2}|w|)^{n} \\
& \leq  \sum_{n=1}^\infty
\frac{1}{(n!)^{\frac{1}{2}}} ||p||_{k+1}^{n} n^{k} (\sqrt{2}|w|)^{n}.
\endalign 
$$ This implies (1) by the definition of $||.||_k$. \par
Ad (2):  By the observation
on Page 267 of [PS] there exists an element $q\in G$ such that
$$
q \exp(v(p)) q^{-1} = \exp(X(p))
$$ up 
to a scalar as operators on $F$.
Since the action of $q$ commutes with the action of rotation, (2)
follows from (1). \par
Ad (3): It is enough to consider the case when $G$ is simple. Since any
level $L$ vacuum representation appears as a direct summand of 
$F^{\otimes L}$, we just have to prove (3) for the representation
$\pi^{\otimes L}$, but this follows immediately from (2). 

\enddemo
\hfill \qed

\proclaim{Lemma 5}
(1)

If $f,p=-p^*$ are smooth functions, then
$$
\pi^{0}(\exp(X(p)))\pi^{0}( V(\psi,f)) x = \pi^{0}(V(\psi,f))
\pi^{0}( \exp(X(p)))  x
$$ for any $x\in H^{0}_0, X\in Lie (H)$;\par
(2) Let $f,p=-p^*$ be smooth functions on $S^1$ with support $f\subset I$
and support $p\subset I'$.  If $X=X_\alpha$  or $X= Y_\alpha$,
then:   
$$
\pi^{0}(V(\psi,f))\pi^{0}( \exp(X(p)))x= 
\pi^{0}( \exp(X(p))) \pi^{0}(V(\psi,f)) x
$$ for any $x\in H^{0}_0$.
\endproclaim
\demo{Proof}
Ad (1): 
If $X\in Lie(H)$, $X(i) \psi = 0$ for any $i\geq 0$ by the
definition of $\psi$,  it follows from (1) of lemma 3
that
$$
[X(p), V(\psi,f)]x = 0
$$ for any $x\in H^0_\infty$.    
By (3) of lemma 4, the condition of (2) of lemma 3 
(note in (2) of lemma 3 $tX(p)= X(tp)$ by definition)
is satisfied,
and the identity follows by using (2) of lemma 3
and the fact that $H^0_0$ is norm dense in $H^0$.\par

Ad (2): Since the support of $f$ and the support of $p$ are 
disjoint, by (1) of lemma 3 
$$
[X(p), V(\psi,f)]x = 0
$$ for any $x\in H^0_\infty$. 
By (3) of lemma 4, the condition of (2) of lemma 3 is satisfied,
and the identity follows by using  (2) of lemma 3 and 
the fact that $H^0_0$ is norm dense in $H^0$. \par

\enddemo
\hfill \qed
\proclaim{Lemma 6}
Let $S$ be the set which consists of elements $\pi^0(\exp(X(p)))$ with
$p=-p^*$ smooth if $X\in Lie(H)$, and $p=-p^*$ smooth,  
support $p \subset I'$ 
if $X=X_\alpha$  or $X=Y_\alpha$. Then the $C^*$ algebra generated by
$S$ is strongly dense in $\pi^0(LH)''\vee \pi^0(L_{I'}G)''$ if $H$
is simply connected.
\endproclaim
\demo{Proof}
Note $S=S^*$. Since every element of $L_{I'}G$ (resp. $LH$) is a product of
exponentials in $L_{I'} g$ (resp. $Lh$), cf. P. 487  of [W2] (we use the fact
that $G,H$ are simply connected here), we just have to show every
element of the form $\pi^0(\exp(X(p)))$ with  $p=-p^*$ smooth,  support
$p \subset I'$,
and $X\in g$ is in the von Neumann algebra $M$ generated by $S$.  
Assume $X=\sum_i c_i X_i$, where $c_i\in {\Bbb R}$ and $X_i$ is either
$X_\alpha$ or $Y_\alpha$. Note that $\pi^0(X(p))$ and $\pi^0(X_i(p))$
are essentially skew self-adjoint operator with a common core $H^0_0$.
By abuse of notations, we will use the same symbol to denote its
closure. 
Let $a\in M'$. Then we have:
$$
\pi^0(X_i(p))a x = a\pi^0(X_i(p)) x
$$ for any $x\in H^0_0$, so
$$
\pi^0(X(p))a x = a\pi^0(X(p)) x
$$ for any $x\in H^0_0$, and it follows that the closure of 
$$
\pi^0(X(p))
$$ is affiliated with $M$, so $\pi^0(\exp(X(p)))$ is in $M$. 
\enddemo
\hfill \qed
\par
\proclaim{Proposition 2.3}
Suppose $H\subset G$, $H$ is  simply connected,  
and $\pi^0$
is the vacuum representation of $LG$.
Let $f$ be a  smooth function with support $f\subset I$. Then
$\pi^0(V(\psi,f))$ is affiliated with von Neumann algebra
$\pi^0(LH)'\cap \pi^0(L_IG)''$.
\endproclaim
\demo{Proof}
By lemma 5, 
$$
\pi^0(V(\psi,f)) u x = u\pi^0(V(\psi,f)) x 
$$ for any $u\in S$ where $S$ is the generating set as in lemma 6 and
$x\in H^0_0$. Note by Haag duality in Prop. 2.1 
$$
\pi^0(L_IG)' = \pi^0(L_{I'}G)''. 
$$
The proposition now follows from (3) of lemma 2 and lemma 6.  
\enddemo  
\hfill \qed
\par
Now we can finish the proof of Th. 2.3.
\demo{Proof of Th. 2.3}
Let $\psi \in
H_{0,0}\otimes \Omega_0$  be an eigenvector of $D$
with eigenvalue $n \in {\Bbb N}$.
Note $\psi = V(-n)\Omega=V(\psi, p) \Omega$ with $p=z^{-1}$. Choose two smooth
functions $f_1$ and $f_2$,  with support $f_1 \subset I_1\in {\Cal I}$ and 
support $f_2 \subset I_2\in {\Cal I}$, 
and $f_1 + f_2 =1$. Then
$$
\psi = V(\psi, p) \Omega = V(\psi, pf_1) \Omega + V(\psi, pf_2) \Omega 
.$$
By Prop. 2.3, the closed operator $V(\psi, f)$ is affiliated with
$\pi^0(LH)'\cap \pi^0(L_JG)''$  if the support $f \subset J \in {\Cal I}$. 
Let $V(\psi, f) = U|V|$ be the polar decomposition.  By lemma 4.4.1
of [Mv], 
$U$ and  $\exp(it|V|), \forall t\in {\Bbb R}$ are in
$\pi^0(LH)'\cap \pi^0(L_JG)''$.  
Since $\Omega$ is in the domain of $|V|$, it
follows by Stone's theorem (cf. P. 266 of [RS]) that
as $t\rightarrow 0$, $(e^{it |V|} \Omega  -\Omega)/t \rightarrow
|V| \Omega$, and this shows
$$
V(\psi, f) \Omega = U|V| \Omega \in 
\overline{\pi^0(LH)'\cap \pi^0(L_JG)''\Omega} 
\subset{\Cal H}
.$$  
So 
$$
V(\psi, pf_1) \Omega \in {\Cal H},  V(\psi, pf_2) \Omega \in {\Cal H}
$$
and  it follows that $\psi \in {\Cal H}$ by the expression given at the
beginning of the proof. 
\enddemo
\hfill \qed  
\par
\subheading {2.3 Two Conjectures}
We will use the notations in 2.1. 
\proclaim{Conjecture 1}
The covariant representations $\pi_{i,\alpha}$ can
be decomposed into a direct sum of a finite number of irreducible
representations
and the localized sectors corresponding to these
irreducible
representations  generate a finite dimensional
ring over  $\Bbb {C}$ under the  product of sectors (cf. 4.1).
\endproclaim
Conjecture 1 comes from the physicists' argument
that the coset $H\subset G$ CFT is a rational CFT: there are only a finite
number of primary fields. Here the primary fields correspond to the
representations or sectors.  In some cases, e.g., $G\subset G\times G$  where
the inclusion is diagonal, there are
also conjectures on the structure constants of the ring (cf. [FKW] and
[BBSS]). More precisely the conjectures in \S3 and 4 of [FKW] are 
about certain
representations
of W-algebras with critical parameters. The W-algebras defined in [FKW]
are closed related to the coset $G \subset 
G_1\times G_m$ (cf:  
[Watts]),
for an example, the representations of  W-algebras in 
[FKW] have the same characters as
those which
come from the coset. We shall call the irreducible conformal 
precosheaves  of the cosets
$SU(N)\subset SU(N)_1 \times SU(N)_m$ {\it coset  $W_N$-algebras with
critical parameters}.   Note that coset $W_2$ algebras with
critical parameters 
are the  irreducible conformal precosheaves 
corresponding to  Virasoro algebras studied in  [GKO]
and [Luke].  
\par
To state conjecture 2, let 
$L_0^{g,h}$ be the generator of rotation group for the coset as 
in the proof of Prop. 2.2.  Then $e^{-\beta L_0}, \beta>0$
is a trace-class operator on $H_{i,\alpha}$ by Th. B of [KW].  Denote by
$d_{i,\alpha}$ the statistical dimension (cf.  
4.1) of $\pi_{i,\alpha}$. Then we
have: (cf. [L5] and (4.27) of [FG])
\proclaim{Conjecture 2 (also known as Kac-Wakimoto formula in [L5])}
$$
d_{i,\alpha} = lim_{\beta \rightarrow 0} \frac{Tr_{H_{i,\alpha}} e^{-\beta L_0
}} {Tr_{H_{0,0}} e^{-\beta L_0}}.      
$$
\endproclaim
Both of these conjectures are highly nontrivial.
The  results in [W2] prove these conjectures
in the case $G$ is of type $A$ and $H$ is a trivial group. For the case
of  coset $W_2$ algebras with
critical parameters the above conjectures follow
from the results of [Luke]. Note  that Conjecture 2
immediately implies Kac-Wakimoto conjecture (cf. Conj. 2.5 in [KW]).
In fact, Conjecture 2 can also be stated as:
$$
d_{i,\alpha} = \frac{b(i,\alpha)}{b(0,0)}
$$ 
where $b(i,\alpha)$ is defined as in \S 2 of [KW] with our 
$(i,\alpha)$ identified with $(\Lambda,\lambda)$ in (2.5.4) of [KW]. 
Conj. 2.5 in
[KW] states that 
$b(i,\alpha) >0 .$ 
Conjecture 2 is stronger than this since $d_{i,\alpha}\geq 1$ 
and
$b(0,0)>0$ by definitions. \par

Note that
the Kac-Wakimoto hypothesis (cf. Page 161 of  [KW]) also implies
Kac-Wakimoto conjecture, but
the first counter-example to Kac-Wakimoto hypothesis has been found in [X2]
by considering subfactors associated with conformal inclusions. So far
Conjecture 2 and hence the Kac-Wakimoto conjecture have been checked to
be true in all known examples.                 

\heading \S 3. Commuting Squares \endheading
We will use the notations of 2.1. All the cosets considered in this
section are assumed to verify the assumptions of Th. 2.3 unless
stated otherwise. 
For the definitions and properties 
of statistical dimensions and minimal index, see
4.1. 
\proclaim{Definition (cofiniteness)}
The coset   $H\subset G_L$ is called cofinite if the inclusion
$$
(\pi^0 ( L_{I}G)'' \cap\pi^0 (L_IH)') \vee \pi^0 (L_IH)'' \subset
\pi^0 ( L_{I}G)''
$$  has finite statistical dimension. The statistical dimension
of the inclusion is denoted by $d(G/H)$.
\endproclaim
Note that $d(G/H)$ does not depend on the choice of $I$ by the
covariance property of representations (cf. Prop. 2.1 of [GL]), and 
we can replace $\pi^0$ by any level $L$ representation of $LG$ in
the above definition due to the local equivalence of these 
representations (cf. Th. II. B of [W2]). \par
Let $\pi^i$ be an irreducible
projective representations of $LG$ with positive energy at level $L$
on Hilbert
space $H^i$.  Recall (cf. 2.1)  when restricting to $LH$, $H^i$ decomposes as:
$$
H^i = \sum_\alpha H_{i,\alpha} \otimes H_\alpha
,$$ and  $\pi_\alpha$ are irreducible projective representations of $LH$ on
Hilbert space $H_\alpha$, 
and the sum is over $\alpha$ such that $(i,\alpha)\in exp$.
Consider the following inclusions:
$$
\align
&(\pi^i ( L_{I}G)'' \cap\pi^i (L_IH)') \vee \pi^i (L_IH)'' \subset
\pi^i ( L_{I}G)'' \subset \pi^i ( L_{I'}G)' \\
&\subset ((\pi^i ( L_{I'}G)'' \cap\pi^i (L_{I'}H)') \vee \pi^i (L_{I'}H)'')'
\endalign
$$
Note that
$$ \pi^i ( L_{I'}G)' \subset
((\pi^i ( L_{I'}G)'' \cap\pi^i (L_{I'}H)') \vee \pi^i (L_{I'}H)'')'
$$
has the same statistical dimension as 
$$
(\pi^i ( L_{I'}G)'' \cap\pi^i (L_{I'}H)') \vee \pi^i (L_{I'}H)'' \subset
\pi^i ( L_{I'}G)''
$$
which is $d(G/H)$. By the multiplicativity of statistical
dimensions
(cf. 4.1) 
the statistical dimension of the inclusion
$$
\align
&(\pi^i ( L_{I}G)'' \cap\pi^i (L_IH)') \vee \pi^i (L_IH)'' \subset
\\
& ((\pi^i ( L_{I'}G)'' \cap\pi^i (L_{I'}H)') \vee \pi^i (L_{I'}H)'')'
\endalign
$$
is $d_id(G/H)^2$,  where $d_i$ is the statistical dimension of
$\pi^i ( L_{I}G)'' \subset \pi^i ( L_{I'}G)'$. 
On the other hand by the additivity of 
 statistical dimension (cf. 4.1)) the statistical dimension of
the above inclusion is $\sum_\alpha d_{(i,\alpha)} d_\alpha$, where
$d_{(i,\alpha)}$ and $d_\alpha$ are the statistical dimensions of
$\pi_{i,\alpha}$ and $ \pi_\alpha$ respectively. So we have
$$
d_i d(G/H)^2 =\sum_\alpha d_{(i,\alpha)} d_\alpha. \tag 3.1
$$ 
When any of the statistical dimensions in formula (3.1) are
$\infty$, then (3.1) is  understood as  the statement that both
sides of the equation  are $\infty$. See the paragraph before Prop. 4.2
for a slightly different derivation of formula (3.1). \par
When $i=0$ is the vacuum representation, the statistical dimension
$d_0$ of the inclusion
$$
\pi^0 ( L_{I}G)'' \subset \pi^0 ( L_{I'}G)'
$$
is $1$ by Haag duality in Prop. 2.1. 
It follows from formula (3.1) that $H\subset G_L$ is cofinite if and only if
$d_{(0,\alpha)} d_\alpha <\infty$ for all
$(0,\alpha)\in exp$.
Hence  Conjecture 2  implies the cofiniteness for any
coset.\par 
For simplicity we will drop the subscript $L$ in the following 
when no confusion 
arises. 
Note that if $H=\{e \}$ is the trivial group , then $H\subset G$ is
cofinite. So the statement 
``if $H_1\subset G$ is cofinite and $H_1\subset H_2\subset G$,
then $ H_2\subset G$ is cofinite"  is as difficult to prove as
the statement ``$ H_2\subset G$ is cofinite" by simply taking
$H_1=\{e \}$. 
But we have:
\proclaim{Proposition 3.1}
Suppose  $H_1\subset H_2\subset G$. \par
(1) If  $H_1\subset G$ is cofinite, then
$H_1\subset H_2$ is cofinite, and $d(G/H_1) \geq d(H_2/H_1)$; \par
(2) If  $H_1\subset H_2$ and $H_2 \subset G$ are cofinite, then
$H_1\subset G$ is cofinite and $d(G/H_1) \leq d(G/H_2) \times 
d(H_2/H_1)$.
\endproclaim
This proposition  is proved below by using commuting squares.
Commuting squares where all the algebras
are finite type can be found in reference [We],[Po]. But
we will consider the case where all the algebras are type $III$. \par

Since
the action of the modular group of 
$\pi^0(L_IG)''$ with respect to the vacuum vector
$\Omega$ is  geometric  and ergodic (cf. 2.1), 
it follows from Takesaki's theorem
(cf. [MT] or P. 495 of [W2]) that  the von Neumann algebras 
$\pi^0(L_IH_i)'' \vee (\pi^0(L_IH_i)'\cap \pi^0(L_IG)''), \ i=1,2$ and
$\pi^0(L_IH_1)'' \vee (\pi^0(L_IH_1)' \cap \pi^0(L_IH_2)'')
 \vee (\pi^0(L_IH_2)'\cap \pi^0(L_IG)'')$ are factors, and 
there exist normal faithful conditional expectations $\epsilon_i:
\pi^0(L_IG)''\rightarrow \pi^0(L_IH_i)'' 
\vee (\pi^0(L_IH_i)'\cap \pi^0(L_IG)'')$
,$\ i=1,2$ and $\epsilon:  
\pi^0(L_IG)''\rightarrow 
\pi^0(L_IH_1)'' \vee (\pi^0(L_IH_1)' \cap \pi^0(L_IH_2)'')
 \vee (\pi^0(L_IH_2)'\cap \pi^0(L_IG)'')$.  Moreover, these conditional
expectations preserve the state $\omega$ on $\pi^0(L_IG)''$ defined by
$\omega(x) = (x \Omega, \Omega), $ i.e.,
$\omega(x) = \omega(\epsilon'(x)),\forall x\in \pi^0(L_IG)''$,
when $\epsilon' = \epsilon_1, \epsilon_2, \epsilon$ respectively.
Then we have:
\proclaim{Lemma 3.1 (Commuting Square)}
(1) $\epsilon_1 \cdot \epsilon_2 =  \epsilon_1 \cdot \epsilon_2 = \epsilon$;
\par
(2) $ \epsilon_1, \epsilon_2$  are minimal if
$\epsilon$ has finite index.
\endproclaim
\demo{Proof}
Assume the vacuum representation $\pi^0$ of $LG$ decomposes with
respect to $LH_2$ as:
$$
H^0 = \oplus_\alpha H_{G/H_2,0,\alpha} \otimes H_{H_2,\alpha}
,$$ then with respect to $LH_1$ the decomposition is:
$$
\align
H^0 = & \oplus_\alpha H_{G/H_2,0,\alpha} \otimes H_{H_2,\alpha} \\
= &  \oplus_{\alpha} H_{G/H_2,0,\alpha} \otimes 
(\oplus_\beta H_{H_2/H_1,\alpha,\beta}
\otimes H_{H_1,\beta})
\endalign
$$ Let  $P_1$ be the projection from $H^0$ onto 
$$ 
\oplus_\alpha H_{G/H_2,0,\alpha} \otimes H_{H_2/H_1,\alpha,0}
\otimes H_{H_1,0},
$$
$P_2$  the projection from $H^0$ onto 
$$ 
H_{G/H_2,0,0} \otimes \oplus_\beta (H_{H_2/H_1,0,\beta}
\otimes H_{H_1,\beta}), 
$$ and 
$P$  the projection from $H^0$ onto 
$$ H_{G/H_2,0,0} \otimes H_{H_2/H_1,0,0}
\otimes H_{H_1,0}.
$$ It follows from definitions 
that $P_1P_2 =P_2 P_1 =P$.  By Th. 2.3 and Reeh-Schlieder
Theorem in Prop. 2.1
$$
P_iH^0 = \overline{\pi^0(L_IH_i)'' \vee (\pi^0(L_IH_i)'\cap 
\pi^0(L_IG)'')\Omega},i=1,2, 
$$ 
and 
$$
PH^0 = \overline{\pi^0(L_IH_1)'' 
\vee (\pi^0(L_IH_1)' \cap \pi^0(L_IH_2)'')
 \vee (\pi^0(L_IH_2)'\cap \pi^0(L_IG)'')\Omega}
.$$ 
So for any $x\in \pi^0(L_IG)''$, we
have:
$$
\epsilon_1 ( \epsilon_2 (x)) \Omega  = P_1 \epsilon_2 (x) \Omega
= P_1 P_2 x \Omega =  P x\Omega = \epsilon (x) \Omega,
$$ 
and similarly
$$
\epsilon_2 ( \epsilon_1 (x)) \Omega = \epsilon (x) \Omega.
$$
Since $\Omega$ is separating for $\pi^0(L_IG)''$, 
(1) of lemma is proved. \par
Note by the remark at the end of 2.2, the inclusions
$$
\pi^0(L_IH_i)'' \vee (\pi^0(L_IH_i)' \cap \pi^0(L_IG)'') \subset
\pi^0(L_IG)''
$$
are irreducible, so by Prop. 4.3 of [L4] 
$\epsilon_i$ are unique and must be the minimal
conditional expectations, $i=1,2$ if the index of
$\epsilon$ is finite.  
\enddemo
\hfill \qed
\demo{Proof of Prop. 3.1}
(1) As in the proof of lemma 3.1, suppose 
 the vacuum representation $\pi^0$ of $LG$ decomposes with
respect to $LH_2$ as:
$$
H^0 = \oplus_\alpha H_{G/H_2,0,\alpha} \otimes H_{H_2,\alpha}
,$$ and let $P_0$ (resp. $P_{00}$)  be the projection from $H^0$ onto 
$H_{G/H_2,0,0} \otimes H_{H_2,0}$ (resp. $H_{G/H_2,0,0} \otimes \Omega_0$
where $\Omega_0$ is the vacuum vector in $ H_{H_2,0}$ ). 
Then 
$$
\align
\pi^0(L_IH_2)'' \vee (\pi^0(L_IH_2)'\cap \pi^0(L_IG)'') & \simeq 
\pi^0(L_IH_2)'' \vee (\pi^0 (L_IH_2)'\cap \pi^0(L_IG)'') P_{0} \\
&\simeq  \pi_0(L_IH_2)'' \otimes (\pi^0 (L_IH_2)'\cap \pi^0(L_IG)'')
P_{00} \\
&\simeq  \pi^0(L_IH_2)'' \otimes (\pi^0 (L_IH_2)'\cap \pi^0(L_IG)'')
\endalign
$$ where $\otimes$ is the tensor product of von Neumann algebras and
$A\simeq B$ means $A$ and $B$ are *-isomorphic, 
since all the algebras above are factors. 
Note the $*$-simorphism above from
$$
\pi^0(L_IH_2)'' \vee (\pi^0(L_IH_2)'\cap \pi^0(L_IG)'')
$$
to  $\pi^0(L_IH_2)'' \otimes (\pi^0 (L_IH_2)'\cap \pi^0(L_IG)'')$ 
maps 
$\pi^0(x)$ to $\pi^0(x) \otimes 1, \forall x\in L_IH_2,$ and
$y$ to $1\otimes y, \forall y\in \pi^0(L_IH_2)'\cap \pi^0(L_IG)''$.  So this
$*$-simorphism maps 
$$
\pi^0(L_IH_1)'' \vee (\pi^0(L_IH_1)' \cap \pi^0 
 (L_IH_2)'') \vee (\pi^0(L_IH_2)'\cap \pi^0(L_IG)'')  
$$
onto 
$$
(\pi^0(L_IH_1)'' \vee (\pi^0(L_IH_1)' \cap \pi^0 
 (L_IH_2)'')) \otimes (\pi^0(L_IH_2)'\cap \pi^0(L_IG)'').
$$ 
It follows that 
the inclusion 
$$
\align
&\pi^0(L_IH_1)'' \vee (\pi^0(L_IH_1)' \cap \pi^0 
 (L_IH_2)'') \vee (\pi^0(L_IH_2)'\cap \pi^0(L_IG)'')   \subset\\
& \pi^0(L_IH_2)'' \vee 
(\pi^0 (L_IH_2)'\cap \pi^0 (L_IG)'')
\endalign
$$ 
is conjugate to
$$
\align
&(\pi^0(L_IH_1)'' \vee (\pi^0(L_IH_1)' \cap \pi^0 
 (L_IH_2)'')) \otimes (\pi^0(L_IH_2)'\cap \pi^0(L_IG)'')   \subset\\
& \pi^0(L_IH_2)'' \otimes 
(\pi^0 (L_IH_2)'\cap \pi^0 (L_IG)''),
\endalign
$$ hence it is irreducible and by 
Cor. 2.2 of  [L6], its  statistical dimension is
 $d(H_2/H_1)$.\par 
By lemma 3.1, the minimal normal faithful conditional expectation
$\epsilon_1$  restricts  to a  normal 
faithful conditional expectation $\eta$
from  
$$  \pi^0(L_IH_2)'' \vee 
(\pi^0 (L_IH_2)'\cap \pi^0 (L_IG)'')
$$ to  
$$
\pi^0(L_IH_1)'' \vee (\pi^0(L_IH_1)' \cap \pi^0 
(L_IH_2)'') \vee (\pi^0(L_IH_2)'\cap \pi^0(L_IG)''). 
$$ 
By Prop. 4.3 of [L4],  $\eta$ is also minimal, and so
$d_{\eta}=d(H_2/H_1) \leq d_{\epsilon_1}=d(G/H_1)$ 
by the definition
of statistical dimension (cf. 4.1).\par
To prove (2), note that $\eta\epsilon_1=\epsilon$ is a minimal conditional
expextation by Cor. 2.2 of [L6], and $d_\epsilon\geq d_{\epsilon_1}$. 
so we have $d(H_2/H_1) d(G/H_2)= d_\epsilon\geq d_{\epsilon_1} =
d(G/H_1),$ where we have used multiplicativity of statistical
dimensions (cf. 4.1).
\enddemo
\hfill \qed
\par
We consider some examples when Prop. 3.1 can be applied.\par
The conformal inclusion $SU(n)_m \times SU(m)_n \subset SU(nm)_1$ has
been considered in [X2] and the decomposition is given in Th. 1 of [ABI]. Let
$H=SU(n)$ be the first factor in the above inclusion. 
\proclaim{Lemma 3.2}
$\pi^0(L_I SU(nm))''\cap \pi^0(L_ISU(n))' = \pi^0(L_I SU(m))''$.
So the irreducible conformal precosheaf of coset $SU(n) \subset SU(nm)_1$ is 
the  irreducible conformal precosheaf of $LSU(m)$ at level $n$.
\endproclaim
\demo{Proof}
From the definition we have:
$$
\pi^0(L_I SU(nm))''\cap \pi^0(L_ISU(n))' \supset \pi^0(L_I SU(m))''
.$$
Since (cf. remark after Prop. 2.2) the action of modular group of 
$$
 \pi^0(L_I SU(nm))''\cap \pi^0(L_ISU(n))'
$$ with respect to the vacuum vector $\Omega$ is geometric and fixes globally
$\pi^0(L_I SU(m))''$, 
by Takesaki's theorem (cf. [MT] or P. 495 of [W2]), we just have to show that
$$
\overline{\pi^0(L_I SU(nm))''\cap \pi^0(L_ISU(n))' \Omega} \subset 
\overline{\pi^0(L_I SU(m))''\Omega}.
$$  
By  the decomposition of 
$H^0$ with respect to $LSU(n) \times LSU(m)$ in Th. 1 of [ABI],
$\Omega = \Omega_{0,0} 
\otimes \Omega_0 \in H_{0,0}\otimes H_{0},$ where 
$H_{0,0}$  
and $ H_{0}$   are vacuum representations 
 of $LSU(m)$ and
$LSU(n)$ respectively, and $ \Omega_{0,0}$, $ \Omega_0$ are vacuum
vectors for $LSU(m)$ and  $LSU(n)$ respectively.
By Reeh-Schlieder
Theorem in Prop. 2.1, 
$\overline{\pi^0(L_I SU(m))''\Omega} =H_{0,0} \otimes \Omega_0$, 
but by the observation before Prop. 2.2 we have
$$
\overline{\pi^0(L_I SU(nm))''\cap \pi^0(L_ISU(n))' \Omega} \subset 
H_{0,0}\otimes \Omega_0 . 
$$ It follows that 
$$
\overline{\pi^0(L_I SU(nm))''\cap \pi^0(L_ISU(n))' \Omega} \subset 
\overline{\pi^0(L_I SU(m))''\Omega},
$$ and the lemma is proved.
\enddemo
\hfill \qed
\par
By lemma 3.2, Th. 1.2 of [X1]  and 
formula (3.1),  
the inclusion $SU(n)_{k+l} \subset SU(n(k+l))_1$ is
cofinite. Since $SU(n)_{k+l} \subset SU(n)_k \times SU(n)_l \subset
SU(n(k+l))_1$ where the first inclusion is diagonal, by (1) of Prop. 3.1
the diagonal inclusion $SU(n)_{k+l} \subset SU(n)_k \times SU(n)_l$
is cofinite, and use (2) of Prop. 3.1 repeatedly we conclude that
the diagonal inclusion $SU(n)_k \subset SU(n)_1 \times...\times
SU(n)_1$ is also cofinite, where there are $k$ factors in the product.
It follows by (1) Prop. 3.1 that 
$SU(n)_{k_1+...+k_m} \subset SU(n)_{k_1} \times...\times SU(n)_{k_m}
$ is cofinite, $k_i\in {\Bbb N}, i=1,...,m$, since
$$
SU(n)_{k_1+...+k_m} \subset SU(n)_{k_1} \times...\times SU(n)_{k_m}
\subset SU(n)_1\times...SU(n)_1
$$ where there are $k_1+...+k_m$ factors in the last group.
\par
Suppose $H_k \subset G_1$ is a conformal inclusion,
 $H$ is simple and of type $A$, $G$ is
simple
 and $k$ is the Dynkin index (cf. P. 170 of [KW]). 
An infinite list can be found in [X2].
Let $l\in {\Bbb N}$. Since $H_{kl} \subset H_k \times...\times H_k$ is
cofinite by the previous paragraph and $H_k \times...\times H_k \subset
G_1\times...\times G_1$ is cofinite by Prop. 2.4 of [X1], it follows by
(2) of
Prop. 3.1 that $H_{kl} \subset G_1\times...\times G_1$ is cofinite, and
by (1) of Prop. 3.1 $H_{kl} \subset G_l$ is cofinite. \par
Finally let us consider the case $H\subset G_m$ with $G=SU(l)$ and
$H$ is the Cartan subalgebra of $G$, a $l-1$ dimensional torus.
We will first consider the inclusion 
$$
H\subset G_1\times...\times G_1
$$
where there are $m$ factors in the product, and the inclusion
is diagonal. Define $\tilde G:= G\times G...\times G$ where
there are $m$ factors in the product. \par
The  irreducible projective representations of
$LH$ at level $m$ have been classified in Prop. 9.5.10 of [PS]. 
Let us describe this result in our case. These irreducible  projective 
representations are in fact representations of ${\Cal L}H$, which is 
a central extension of $LH$ induced from the central extension 
${\Cal L}G$ of 
$LG$ (cf. P. 483 of [W2] or Chap. 4 of [PS]). 
Write $LH\simeq \text{\rm Hom}(S^1, H) \times H\times V$, where $H$ is the
subgroup of constant loops, and $V$ is the vector space of maps
$f:S^1\rightarrow Lie(H)$ with integral $0$, which is regarded as 
a subgroup of $LH$ by the exponential map. The 
dentity component of ${\Cal L}H$ is canonically a product $H\times \tilde V$,
where $\tilde V$ is the Heisenberg group defined by a skew form
on $V$.
The center of
the identity component of ${\Cal L}H$ is $H\times S^1$.  
Let 
$\xi=(\xi_1,..., \xi_{l-1},
(\xi_1...\xi_{l-1})^{-1}) \in LH, $ 
where each $\xi_i\in C^\infty(S^1, S^1)$ 
has winding number $x_i, i=1,...,l-1$. 
The conjugate action of $\xi$ 
on the center 
$H\times S^1$ of the identity component of ${\Cal L}H$ is given by
(cf. P. 192 of [PS])
$$
(t,u)\rightarrow (t, u t_1^{x_1+(x_1+...x_{l-1})}... 
t_{l-1}^{x_{l-1}+(x_1+...x_{l-1})})
$$ where $t=(t_1,...,t_{l-1}, (t_1...t_{l-1})^{-1})\in H$. 
Introduce an equivalent relation  on ${\Bbb Z}^{l-1}$ by:
$ (n_1,...,n_{l-1}) \sim  (n_1',...,n_{l-1}')$ iff
there exists $(m_1,...m_{l-1}) \in {\Bbb Z}^{l-1}$ with
$m_1+...+m_{l-1} \in l{\Bbb Z}$ such that
$(n_1',...,n_{l-1}') = (n_1+mm_1,...n_{l-1}+mm_{l-1})$. Denote
the equivalence class of $(n_1,...,n_{l-1})$ by 
$
[n_1,...,n_{l-1}]
$
or simply $[n]$. The  irreducible  representation of
${\Cal L}H$ at level $m$ on Hilbert space $H_{[n]}$ has the following form
(cf. P. 192 of [PS]):
$$
H_{[n]}= \oplus_{(a_1,...,a_{l-1})\sim (n_1,...,n_{l-1})} 
H_{(a_1,...,a_{l-1})},
$$ where on $H_{(a_1,...,a_{l-1})}$, the center $H\times S^1$ of
the identity componenet of ${\Cal L}H$ acts as
$(t,u)\rightarrow t_1^{a_1}... t_{l-1}^{a_{l-1}}u \times \text{\rm id}$, 
and on
$H_{(a_1,...,a_{l-1})}$ the representation of the Heisenberg group $\tilde V$
is irreducible (and unique by Prop. 9.5.10 of [PS]).
\par
When restricting to ${\Cal L}H$, 
the vacuum representation $\pi^0$ of ${\Cal L}\tilde G:={\Cal L}G\times
... {\Cal L}G$ (there are m factors in the product) on 
$(H_v)^{\otimes m}$ decomposes  as:
$$
(H_v)^{\otimes m} = \sum_{l|\sum_i n_i} H_{0,[n]} \otimes
H_{[n_1,...n_{l-1}]} 
,$$ and $\pi_{[n]}$ are  irreducible projective representations
of ${\Cal L}H$ on  $H_{[n_1,...n_{l-1}]}$ (cf. \S2.6 of [KW]). 
Let $\alpha \in LH \times...\times LH$ be a loop 
of the form $\xi \times 1\times...\times 1$, where 
$\xi=(\xi_1,..., \xi_{l-1},
(\xi_1...\xi_{l-1})^{-1}) \in LH, $ 
and  each $\xi_i\in C^\infty(S^1, S^1) $ 
has  winding number $x_i, i=1,...,l-1$. 
We can assume that
$\alpha$ is localized on $I$. Define 
$$
Ad_\alpha.y:= \alpha y\alpha^{-1},
Ad_\alpha.\pi^0 (y) = \pi^0 (\alpha) \pi^0 (y) \pi^0 (\alpha^{-1}),  
\forall y\in {\Cal L}\tilde G.
$$ 
Note (remember that ${\Cal L}H$ is diagonally included in ${\Cal L}\tilde G$)
$$
Ad_\alpha.y \in {\Cal L}H , Ad_\alpha.\pi^0 (y) \in \pi^0({\Cal L}H)
, \forall y\in {\Cal L}H.
$$ 
Then by definitions we have:
$$
\pi_{[n]} (Ad_\alpha y) \simeq \pi_{[n+b]} (y), \forall y\in {\Cal L}H
,$$ where $[n+b]= [n_1+b_1,...,n_{l-1}+ b_{l-1}]$, with 
$b_i= x_i +(x_1+...x_{l-1}), i=1,...,l-1$.
Note that this implies that $\pi_{[n]}$ has statistical dimension
$1$ since $Ad_\alpha$ is a localized automorphism. 
Also note that since  $Ad_\alpha$ is an automorphism of $\pi^0({\Cal L}H)$,  it
is also an automorphism of $\pi^0({\Cal L}H)'\cap \pi^0({\Cal L}_J\tilde G)''$ 
for any interval $J$. \par
We claim that
$$
\pi_{0,[n]} (Ad_\alpha.y) \simeq 
 \pi_{0,[n+b]} (y)
$$ for any $y\in A(J)$, where $A(J)$ is the conformal precosheaf  for the
coset $H\subset \tilde G$.
In fact let $U_{[n]}: H_{[n]}\rightarrow H_{[n+b]}$ be  a unitary
map intertwinning the action of 
${\Cal L}H$ and the action of $Ad_\alpha. {\Cal L}H$, 
and let $W= V\otimes U: H^0 \rightarrow H^0$
be a unitary map such that $V\otimes U( z\otimes y)= V_{[n]}z \otimes
U_{[n]}y \in    H_{0,[n+b]} \otimes  H_{[n+b]}$ for any
$z\otimes y \in  H_{0,[n]} \otimes  H_{[n]}$. It follows that
$$
W^* \pi^0(\alpha) \in \pi^0({\Cal L}H)' = \oplus_{[n]} B(H_{0,[n]}) \otimes
id_{H_{[n]}}
,$$ so we have
$$
\pi^0(\alpha) = \oplus_{[n]}V'_{[n]}\otimes U_{[n]}
$$ where $  V'_{[n]}:  H_{0,[n]}\rightarrow H_{0,[n+b]}$ is unitary.
Note
$
\pi^0 (Ad_\alpha.y) = \pi^0(\alpha) \pi^0(y) \pi^0(\alpha)^{-1},
$ and 
$$
\pi^0(y) \in \oplus_{[n]} B(H_{0,[n]}) \otimes
id_{H_{[n]}}, \forall y\in A(J).
$$
Hence
$$
\pi_{0,[n]} (Ad_\alpha.y) =
  {V'}_{[n]} \pi_{0,[n+b]} (y)  {V'}_{[n]}^*, \forall y\in A(J).
$$ 
Now choose $x$ so that
$b_i= -n_i, i=1,...,l-1$, we get
$$
\pi_{0,[n]} (Ad_\alpha.y) =
  {V'}_{[n]} \pi_{0,[0]} (y)  {V'}_{[n]}^*
.$$  

So $\pi_{0,[n]}$ has the same statistical dimension as $\pi_{0,[0]}$
since $Ad_\alpha$ is a localized automorphism.

\par
If $\pi_{0,[0]}$ is the vacuum representation,
then $\pi_{0,[n]}$ has  statistical dimension 1. Note $\pi_{[n]}$ also
has  statistical dimension 1, by formula (3.1)
we  conclude that the diagonal inclusion 
$H\subset G_1\times...\times G_1$ is cofinite. \par  
We claim that $\pi_{0,[0]}$ is indeed the vacuum representation. Note
this does not follow directly from Th. 2.3 since we assume $H$ is
simply connected in the theorem. However, the assumption that
$H$ is
simply connected is only used in the proof of lemma 6.  
From the proof of lemma 6, we  see that the smeared vertex operators
in Prop. 2.3 are affiliated with  von Neumann algebra
$$
\pi^0((LH)^0)' \cap \pi^0(L_I\tilde G)''
,$$ where $(LH)^0$ is the connected component of $LH$ that contains
identity. Note that $LH$ is generated as a group by $(LH)^0$ and
a set of elements with non-trivial winding numbers, and we can
certainly choose these elements to be in $L_{I'}H \subset L_{I'} \tilde G$.
So
$$
\pi^0(LH)'' \subset  \pi^0((LH)^0)'' \vee \pi^0(L_{I'} \tilde G)''.
$$
Hence if $p\in \pi^0((LH)^0)' \cap \pi^0(L_I \tilde G)''$, then 
$p\in  \pi^0(LH)' \cap \pi^0(L_I \tilde G)''$. On the other hand
 
$$
\pi^0((LH)^0)' \cap \pi^0(L_I \tilde G)'' \supset   \pi^0(LH)' 
\cap \pi^0(L_I \tilde G)'' 
.$$ 
So 
$$
\pi^0((LH)^0)' \cap \pi^0(L_I \tilde G)'' =   \pi^0(LH)' 
\cap \pi^0(L_I \tilde G)'' 
.$$ 

This shows that Prop. 2.3, and therefore Th. 2.3 hold for any 
pair $H\subset \tilde G$ as long as $\tilde G$ is semisimple and simply 
connected.
It follows now that Prop. 3.1 can be applied to the present case 
for $H\subset G_m\subset G_1\times... \times G_1$ since we only
use Th. 2.3 in its proof. So we conclude that $H\subset G_m$ is 
cofinite. \par

To summarize, we have proved the following:
\proclaim{Corollary 3.1}
The following inclusions are cofinite: \par
(1) $G_{k_1+k_2+...+k_m} \subset G_{k_1} \times... \times
G_{k_m}$ where the inclusion is diagonal, $k_i\in {\Bbb N}, i=1,...,m
$ and $G=SU(n)$;\par
(2) $H_{lk} \subset G_l$, if  $H_{k} \subset G_1$
is a conformal inclusion where $k$
is the Dynkin index, $l\in {\Bbb N}$, $H$ is simple and of type $A$ and $G$ is
simple ; \par
(3) $H\subset G_m$, where $H$ is the Cartan subgroup of $G$.
\endproclaim  
\heading{\S4. Braided endomorphisms}\endheading
All the cosets considered in this section are assumed to 
verify the assumptions of Th. 2.3 unless stated otherwise.
\subheading{4.1 Some results from [X1]}
In this subsection we recall some of the results from [X1] which
will be used in 4.2. We start with some preliminaries on
sectors to set up notations. \par
Let $M$ be a properly infinite factor
and  $\text{\rm End}(M)$ the semigroup of
 unit preserving endomorphisms of $M$.  In this paper $M$ will always
be a type $III_1$ factor.
Let $\text{\rm Sect}(M)$ denote the quotient of $\text{\rm End}(M)$ modulo
unitary equivalence in $M$. 
It follows from
\cite{L3} and \cite{L4} that $\text{\rm Sect}(M)$  is endowed
with a natural involution $\theta \rightarrow \bar \theta $, and
$\text{\rm Sect}(M)$ is
a semiring: i.e., there are two operations $+, \times$ on 
$\text{\rm Sect}(M)$ which verifes the usual axioms. The multiplication
of  sectors is simply the composition of sectors. Hence if
$\theta_1, \theta_2$ are
two sectors, we shall write $\theta_1\times \theta_2$  as
$\theta_1\theta_2$.
In [X1],   the image of
$\theta \in \text{\rm End}(M)$ in  $\text{\rm Sect}(M)$ is denoted by 
$[\theta]$. However, since we will be mainly concerned with the ring
structure of certain sectors in section 4, we will denote $[\theta]$
simply by $\theta$ if no confusion arises. 
\par
Assume $\theta \in \text{\rm End}(M)$, and there exists a normal
faithful conditional expectation
$\epsilon:
M\rightarrow \theta (M)$.  We define a number $d_\epsilon$ (possibly
$\infty$) by:
$$
d_\epsilon^{-2} :=\text{\rm Max} \{ \lambda \in [0, +\infty)|
\epsilon (m_+) \geq \lambda m_+, \forall m_+ \in M_+
\}$$ (cf. [PP]).\par         
If $d_\epsilon < \infty$ for some $\epsilon$, we say $\theta$
has finite index or statistical dimension. In this case 
we define
$$
d_\theta = \text{\rm Min}_\epsilon \{ d_\epsilon |  d_\epsilon < \infty \}.
$$   $d_\theta$ is called the {\it statistical dimension} of  $\theta$. 
$d_\theta^2$ is called the {\it minimal index} of $\theta$. In fact in this
case there exists a unique $\epsilon_\theta$ such that
$ d_{\epsilon_\theta} = d_\theta$. $\epsilon_\theta$ is called the
{\it minimal conditional expectation}. 
It is clear
from the definition that  the statistical dimension  of  $\theta$ depends only
on the unitary equivalence classes  of  $\theta$.
When $N\subset M$ with $N\simeq M$, we choose $\theta \in \text{\rm End}(M)$
such that $\theta(M)=N$.   The statistical dimension 
(resp. minimal index) of the inclusion $N\subset M$ is 
defined to be the statistical dimension (resp. minimal index) of $\theta$.\par
Let  $\theta_1, \theta_2\in Sect(M)$. By Th. 5.5 of [L3],  
$d_{\theta_1+\theta_2}= d_{\theta_1} + d_{\theta_2}$, and
by Cor. 2.2 of [L6], $d_{\theta_1\theta_2}= d_{\theta_1}d_{\theta_2}$.
These two properties are usually referred to as the {\it additivity}
and {\it multiplicativity} of statistical dimensions. Also note
by Prop. 4.12 of [L4] $d_\theta = d_{\bar \theta}$. If a sector
does not have finite statistical dimension in any of the
above three equations, then the equation is understood as the 
statement that both sides of the equation are $\infty$. \par  
Assume $\lambda $, $\mu,$ and $ \nu \in \text{\rm End}(M)$ have
finite statistical dimensions.  Let
$\text{\rm Hom}(\lambda , \mu )$ denote the space of intertwiners from
$\lambda $ to $\mu $, i.e. $a\in \text{\rm Hom}(\lambda , \mu )$ iff
$a \lambda (p) = \mu (p) a $ for any $p \in M$.
$\text{\rm Hom}(\lambda , \mu )$  is a finite dimensional vector
space and we use $\langle  \lambda , \mu \rangle$ to denote
the dimension of this space.  Note that $\langle  \lambda , \mu \rangle$
depends
only on $[\lambda ]$ and $[\mu ]$. Moreover we have
$\langle \nu \lambda , \mu \rangle =
\langle \lambda , \bar \nu \mu \rangle $,
$\langle \nu \lambda , \mu \rangle    
= \langle \nu , \mu \bar \lambda \rangle $ which follows from Frobenius
duality (See \cite{L2} or \cite{Y}).  We will also use the following
notation: if $\mu $ is a subsector of $\lambda $, we will write as
$\mu \prec \lambda $  or $\lambda \succ \mu $.  A sector
is said to be irreducible if it has only one subsector. \par 
Let $\theta_i, i=1,...,n$ be a set of irreducible 
sectors with finite index. The
ring generated by $\theta_i, i=1,...,n$ under compositions 
is defined to be a vector space 
(possibly infinite dimensional) over
${\Bbb C}$ with a  basis $\{ \xi_j , j\geq 1\}$, such that $\xi_j$ are
irreducible sectors, $\xi_j\neq \xi_{j'}$ if $j\neq j'$, and the set
$\{\xi_j, j\geq 1\}$ is a list of all   
irreducible sectors which appear as subsectors of 
finite products of $\theta_i,
i=1,...,n$. The ring 
multiplication  on the vector space is obtained naturally from that of
$Sect(M)$.\par
Let $ M(J), J\in {\Cal I}$ 
be an irreducible conformal precosheaf  on Hilbert space $H^0$. 
Suppose
$N(J), J\in {\Cal I}$ is an irreducible conformal precosheaf 
 and $\pi^0$ is a covariant representation
of $N(J)$ on $H^0$ such that $\pi^0(N(J)) \subset M(J)$ is a directed
standard net as defined in Definition 3.1 of [LR] for any directed
set of intervals. Fix an interval $I$ and denote by $N:= N(I), M:=M(I)$.   
For any covariant representation $\pi_\lambda$ (resp. $\pi^i$) of
the irreducible conformal precosheaf 
$N(J),   J\in {\Cal I} $ (resp. 
$M(J), J\in {\Cal I}$), let $\lambda$ (resp. $i$) be the 
corresponding  
endomorphism  of $N$ (resp. $M$) as defined in \S2.1 of [GL]. 
These endomorphisms are obtained by localization in \S2.1 of [GL]
and will be referred to  as  localized 
endomorphisms for convenience. The corresponding sectors will be
called localized sectors.  See the paragraph after the proof 
of Lemma 4.2 for examples.\par
We will use $d_\lambda$ and  $d_i$ to denote the 
statistical dimensions of $\lambda$ and  $i$ respectively. 
$d_\lambda$ and  $d_i$
are also called  the 
statistical dimensions of $\pi_\lambda$ and $\pi^i$ respectively, and they are
independent of the choice of $I$ (cf. Prop. 2.1 of [GL]). \par
Let $\pi^i$
be a covariant representation of  $M(J),  J\in {\Cal I}$
which decomposes as:
$$
\pi^i = \sum_{\lambda} b_{i\lambda} \pi_\lambda
$$ when restricted to  $N(J),  J\in {\Cal I}$, where the sum is finite
and $ b_{i\lambda} \in {\Bbb N}$. Let $\gamma_i:=
\sum_{\lambda} b_{i\lambda} \lambda$ be the corresponding sector of $N$. 
It is shown (cf. (1) of
Prop. 2.8 in [X1]) that there are sectors $\rho, \sigma_i \in
\text{\rm Sect}(N)$ such that:
$$
\rho \sigma_i \bar \rho =\gamma_i.
$$
Notice that  $\sigma_i$ are in one-to-one correspondence with
covariant representations $\pi^i$, and in fact the map $i\rightarrow
\sigma_i$ is an isomorphism of the ring generated by $i$ and
the ring generated by $\sigma_i$.  The subfactor 
$\bar{\rho}(N) \subset N$ is conjugate to $\pi^0(N(I)) \subset M(I)$
(cf. (2) of Prop. 2.6 in [X1]). \par
Now we assume  $\pi^0(N(I)) \subset M(I)$ has finite index. Then for
each  localized sector $\lambda$ of $N$ there exists a sector
denoted by $a_\lambda$ of $N$ such that the following theorem is
true   (cf. [X1]):
\proclaim{Theorem 4.1}
(1)
The map $\lambda \rightarrow a_\lambda$ is a ring homomorphism;\par
(2)
$\rho a_\lambda = \lambda \rho, a_\lambda \bar \rho = \bar \rho \lambda
, d_\lambda = d_{a_\lambda}$;\par
(3) $\langle \rho a_\lambda,  \rho a_\mu \rangle =
\langle a_\lambda,  a_\mu \rangle = \langle  a_\lambda \bar \rho,  a_\mu
\bar \rho  \rangle$; \par
(4)  $\langle \rho a_\lambda, \rho \sigma_i  \rangle =   
\langle a_\lambda, \sigma_i  \rangle =
\langle a_\lambda \bar \rho, \sigma_i \bar \rho
\rangle $ ; \par
(5) (3) (resp. (4)) remains valid if $ a_\lambda,  a_\mu$ (resp.  $a_\lambda$)
is replaced by any of its subsectors; \par
(6) $ a_\lambda \sigma_i =  \sigma_i a_\lambda$.
\endproclaim
\demo{Proof}
(1) to (4) follows from  Th. 3.1, 3.3, Cor. 3.2, lemma 3.4, 3.5 of
 [X1],  (5) is proved on P. 9 of [X2], and (6) is proved on P. 387
of [X1]. It should be noted that these results in \S3  of
[X1] are stated for
conformal inclusions, but  all the
proof there applies verbatim to the present setting.  
\enddemo 
\hfill \qed
  
\subheading{4.2 The  ring structure}
We will apply the results of 4.1 to  the case when $N(I)= A(I)\otimes
\pi_0(L_IH)''$ and $M(I)=\pi^0(L_IG)''$ under the assumption that
$H\subset G_L$ is cofinite, where
$A(I)$ is as in Prop. 2.2 for the coset $H\subset G_L$, and
$\pi_0$ denotes the vacuum representation of $LH$.  
Note that if $H\subset G_L$ is cofinite, 
then $\pi^0(N(I)) \subset M(I)$ has finite index.  
By Th. 4.1, for every localized endomorphisms $\lambda $ of
$N(I)$ we have a map $a: \lambda \rightarrow a_\lambda$ which
verifies (1) to (6) in Th. 4.1. \par
\proclaim{Tensor Notation}
Let $\theta \in End( A(I)\otimes
\pi_0(L_IH)'')$. We will denote $\theta$ by $\rho_1\otimes \rho_2$
if 
$$
\theta(p\otimes 1)= \rho_1(p)\otimes 1, \forall p\in A(I),
\theta(1\otimes p')= 1\otimes \rho_2(p'), \forall p'\in\pi_0(L_IH)'' ,  
,$$ where $\rho_1 \in End( A(I)),\rho_2 \in End(\pi_0(L_IH)'')$.
\endproclaim
\proclaim{Lemma 4.2}
(1) If 
$\theta = \rho_1 \otimes \rho_2$, and 
$$
[\rho_1] =\sum_i [\rho_{1i}], [\rho_2] =\sum_j [\rho_{2j}]
,$$ where all the summations are finite. Then:
$$
 [\theta]= \sum_{i,j}[\rho_{1i} \otimes  \rho_{2j}]
;$$\par
(2)
$$
\langle \rho_1 \otimes \rho_2, \sigma_1 \otimes \sigma_2\rangle
= \langle \rho_1, \rho_2\rangle \langle \sigma_1, \sigma_2\rangle 
,$$ where $\rho_1, \sigma_1$ are in $End(A(I))$, and
 $\rho_2, \sigma_2$ are in $ End(\pi_0(L_IH)'')$.
\endproclaim
\demo{Proof}
(1) follows immediately from the definitions. By (1), we just
have to show (2) in the case that  $\rho_1, \sigma_1,\rho_2, \sigma_2 $
are irreducible sectors. It is obvious that if
$\rho_1 \simeq \sigma_1,\rho_2 \simeq \sigma_2$ as sectors, then 
$\rho_1 \otimes \rho_2 \simeq \sigma_1 \otimes \sigma_2$ as sectors
of $N(I)$. Now suppose 
$\rho_1 \otimes \rho_2 \simeq \sigma_1 \otimes \sigma_2$ as sectors 
of $N(I)$. This means there exists a unitary $u\in N(I)$ such that:
$$
u\rho_1(p) \otimes \rho_2 (p') = \sigma_1(p) \otimes \sigma_2(p')u  
,$$ for any $p\in A(I), p'\in \pi_0(L_IH)''$. By 
the statement on P. 123 of [Stra], there
exists normal  conditional expectation $E: N(I)\rightarrow A(I)\otimes 1$
such that $E(u) \neq 0$. Applying $E$ to the above equation and setting $p'=1$,
we have:
$$
E(u)\rho_1(p)=\sigma_1(p) E(u)
.$$ Since $ \rho_1,\sigma_1$ are irreducible and $E(u)\neq 0$, it follows
that $ \rho_1 \simeq \sigma_1$ as sectors. Similarly one can show that
$ \rho_2 \simeq \sigma_2$ as sectors. This proves (2). 
\enddemo
\hfill \qed
\par

Recall from 2.1  $\pi_{i,\alpha}$ of $A(I)$ are obtained 
in the decompositions of $\pi^i$ of $LG$ with respect to subgroup $LH$,
and we  denote the set of such $(i,\alpha)$ by $exp$. 
For any $J\in {\Cal I}$, 
let $U(J)$ be a unitary operator from $H_{i,\alpha}$ to
$H_{0,0}$ such that:
$$
\pi_{i,\alpha}(p) = U(J)^* \pi_{0,0}(p) U(J), \forall p\in A(J)
.$$ 
Recall $I$ is a fixed interval. Identify $H_{i,\alpha}$ with
$H_{0,0}$ by $U(I')$, we may  choose a representation unitarily
equivalent to $\pi_{i,\alpha}$, still denoted by $\pi_{i,\alpha}$ on
$H_{0,0}$, with the property that
$\pi_{i,\alpha}(p')=p', \forall p'\in A(I')$. It follows that
$\pi_{i,\alpha}(A(I))$ commutes with $ A(I')$.  
By Haag
duality in Prop. 2.1, 
$\pi_{i,\alpha} (p) \in A(I), \forall p\in A(I)$, and so
$
\pi_{i,\alpha}|A(I) \in \text{\rm End}(A(I)).
$ We will denote $\pi_{i,\alpha}|A(I)$ by $(i,\alpha)$. The corresponding
sector in $\text{\rm Sect}(A(I))$ is also 
denoted by $(i,\alpha)$ when no confuison
arises. Note that  $(i,\alpha)$ is an irreducible sector if and only if
$\pi_{i,\alpha}$ is an irreducible covariant representation, since the 
coset conformal precosheaf $A(J), \forall J\in {\Cal I}$ 
is stronly additive by the remarks in 2.1 after
Prop. 2.2. In fact suppose $(i,\alpha)$ is an irreducible sector. Let
$p\in (\vee_{J\in {\Cal I}}\pi_{i,\alpha}(A(J)))'$. Then 
$p\in \pi_{i,\alpha}(A(I'))'= A(I)$. It follows that $p\in \text{\rm Hom} 
((i,\alpha),
(i,\alpha))= {\Bbb C}$ since $(i,\alpha)$ is irreduicble. On the other hand
if $\pi_{i,\alpha}$ is irreducible, and $p\in  \text{\rm Hom} ((i,\alpha),
(i,\alpha))$. Then $p\in A(I)$ and so
$p\in (\pi_{i,\alpha}(A(I')) \vee  \pi_{i,\alpha}(A(I)))'$. But 
$$ 
\pi_{i,\alpha}(A(I')) \vee  \pi_{i,\alpha}(A(I))= 
\vee_{J\in {\Cal I}}\pi_{i,\alpha}(A(J))
$$ by the strong additivity of the coset conformal precosheaf, so
$$p\in (\vee_{J\in {\Cal I}}\pi_{i,\alpha}(A(J)))' = {\Bbb C}
$$ since
$\pi_{i,\alpha}$ is irreducible. Similarly one can show that $
(i,\alpha)\succ (j,\beta)$ if and only of $\pi_{j,\beta}$ appears
as a direct summand of $\pi_{i,\alpha}$,  and $(i,\alpha)$ is equal
to $(j,\beta)$ as sectors if and only  $\pi_{i,\alpha}$ is
unitarily equivalent to  $\pi_{j,\beta}$.\par
Given  $(i,\alpha)\in \text{\rm End}(A(I))$ as above, we define
$(i,\alpha) \otimes 1 \in \text{\rm End} (N(I))$ so that:
$$
(i,\alpha)\otimes 1(p\otimes p') = (i,\alpha)(p)\otimes p' 
, \forall p\in A(I), p'\in \pi_0(L_IH)''
.$$ It is
easy to see that $(i,\alpha)\otimes 1$ corresponds to the
covariant representation $\pi_{i,\alpha} \otimes \pi_0$ of
$N(I)$.
Note that this notation agrees with our tensor notation above. 
Also note that for any covariant representation $\pi_x$
of $A(I)$, we can define a localized sector $x\otimes 1$ of
$N(I)$ in the same way as in the case when $\pi_x=\pi_{i,\alpha}$. 
\par
Each covariant representation $\pi^i$ of $LG$ 
gives rise to an endomorphism
$\sigma_i\in \text{\rm End} (N(I))$ and (cf. subsetion 4.1)
$$
\rho \sigma_i \bar{\rho} = \gamma_i = \sum_{\alpha}
(i,\alpha) \otimes (\alpha)
$$ where the summation is over those $\alpha$ such that $(i,\alpha) \in exp$.
So by the properties 
 of statistical dimensions (cf. 4.1) 
$d_i d_\rho^2= \sum_{\alpha} d_{(i,\alpha)}d_\alpha$. 
Note that this is in fact formula
(3.1), with $d_\rho= d(G/H)$ by definition.  
\proclaim{Proposition 4.2}
Assume $H\subset G_L$ is cofinite.  
We have:\par 
(1)
Let $x, y$ be localized 
sectors of $A(I)$ with finite index. Then
$$
\langle x, y \rangle = 
\langle a_{x\otimes 1}, a_{y \otimes 1} \rangle 
;$$\par
(2) If $(i,\alpha) \in exp$, then
$ a_{(i,\alpha)\otimes 1}  \prec a_{1\otimes \bar{\alpha}} \sigma_i$;  
\par
(3) Denote by $d_{(i,\alpha)}$ the statistical dimension of
$(i,\alpha)$. Then $d_{(i,\alpha)} \leq d_i d_\alpha$, where 
$d_i$ (resp. $d_\alpha$) is the statistical dimension of $i$ (
resp. $\alpha$). 
\endproclaim
\demo{Proof}
Ad (1): 
By the assumption and Th. 4.1, we have
$$
\align
\langle a_{x\otimes 1}, a_{y \otimes 1} \rangle & = 
\langle \rho a_{x\otimes 1} , \rho a_{y \otimes 1} \rangle 
\\
&=\langle (x\otimes 1) \rho,  (y\otimes 1) \rho \rangle \\
&=\langle x\otimes 1 ,  (y \otimes 1) \rho \bar{\rho}\rangle \\ 
&=\langle x\otimes 1 ,  (y \otimes 1) \sum_\delta (0,\delta)
\otimes \delta\rangle \\
&=\langle x\otimes 1 , \sum_\delta y(0,\delta)
\otimes \delta\rangle \\
&= \sum_\delta \langle x ,  y(0,\delta) \rangle
\times \langle 1, \delta \rangle ,
\endalign
$$ where in the last identity we used (2) of lemma 4.2.
Note $\langle 1, \delta \rangle$ is equal to 1 iff $\delta$ corresponds to 
the vacuum representation of $LH$. When $\delta$ is the
vacuum representation, $(0,\delta)$ corresponds to the representation
$\pi_{0,0}$ of $A(I)$, which by Th. 2.3, is the vacuum representation
of $A(I)$, and corresponds to the identity sector. So we have:
$$
\sum_\delta\langle x ,  y(0,\delta) \rangle
\times \langle 1, \delta \rangle = 
\langle x ,  y \rangle
,$$ and the proof of (1) is complete. \par
Ad (2): Since $\pi^i$ has finite statistical dimension and
$H\subset G$ is cofinite, by formula (3.1) $d_{(i,\alpha)} < \infty, 
d_\alpha< \infty, \forall (i,\alpha)\in exp$. So
we can assume $(i,\alpha) = \sum_j m_j x_j$, where the sum
is finite, $m_j\in {\Bbb N},$ and $x_j$ is irreducible and has finite 
index. 
Note that $x_j$ is a localized sector of $A(I)$ (cf. Prop. 2.2
of [GL]), so $a_{x_j\otimes 1}$ is well defined, and it follows from
(1) that $a_{x_j\otimes 1}$ is also irreducible. By Th. 4.1
$$
a_{(i,\alpha)\otimes 1} = \sum_j m_ja_{x_j\otimes 1}  
.$$ By using Th. 4.1 and lemma 4.2 we have:
$$
\align
\langle a_{x_i\otimes 1},  a_{1\otimes \bar \alpha} \sigma_i \rangle &=
\langle a_{x_i\otimes 1}  a_{1\otimes \alpha},  \sigma_i \rangle \\
&= \langle a_{x_i \otimes \alpha},  \sigma_i \rangle \\
&=  \langle a_{x_i\otimes  \alpha} \bar \rho,  
\sigma_i \bar \rho\rangle \\
&=  \langle \bar \rho x_i\otimes  \alpha,  \sigma_i  
\bar \rho \rangle\\
&=  \langle x_i\otimes  \alpha,  \rho \sigma_i \bar \rho\rangle \\
&=  \langle x_i\otimes  \alpha, \sum_\beta (i,\beta) \otimes \beta
\rangle \\
&= \langle x_i, (i,\alpha) \rangle = m_i.
\endalign
$$ 
This shows 
$$
a_{(i,\alpha)\otimes 1} = \sum_j m_ja_{x_j\otimes 1} \prec  
a_{1\otimes \bar{\alpha}} \sigma_i   
.$$ \par
(3) follows immediately from (2) and the fact that $(i,\alpha)$ and
$a_{(i,\alpha)\otimes 1}$ have the same statistical dimension  by 
(2) of Th. 4.1.
\enddemo 
\hfill \qed
\proclaim{Theorem 4.2}
Suppose $H\subset G_L$ is cofinite,  and every irreducible representation
$\pi_\alpha$ of  $LH$ has
finite index, and the localized 
sectors $\{ \alpha \}$ generate a finite
dimensional ring over $\Bbb {C}$ 
under  compositions. Then
Conj. 1 of 2.3 is true.
\endproclaim
\demo{Proof}
By (3) of Prop. 4.2, $d_{(i,\alpha)} < \infty, \forall (i, \alpha)\in exp,$
since $d_i <\infty, d_\alpha <\infty$ by our assumption.
Hence each $(i, \alpha)$ decomposes into a direct sum of a finite number
of irreducible sectors. By Prop. 2.2 of [GL], each $\pi_{i,\alpha}$
decomposes into a direct sum of a finite number
of irreducible covariant representations.\par
By (6) of Th. 4.1, 
$ a_{1\otimes \bar{\alpha}} \sigma_i = \sigma_i a_{1\otimes \bar{\alpha}}
$, so the ring $Y$ generated by   irreducible subsetors of 
$ a_{1\otimes \bar{\alpha}} \sigma_i, \forall
\alpha, \forall i$ has  finite dimension over $\Bbb {C}$
by (1) of Th. 4.1 and the the assumption of the theorem.  
Denote by $X$ the ring generated (under compositions) by the 
set of all irreducible sectors which appear as subsectors
of $(i, \alpha), \forall (i, \alpha)\in exp$. By Prop. 4.2, the 
map $x\in X\rightarrow a_{x\otimes 1}$ is an injective homomorphism 
from $X$ into $Y$. It follows that $X$ is finite dimensional
over $\Bbb {C}$,  and  Conj. 1 of
2.3 is proved.
\enddemo \hfill \qed
\par
By Cor. 3.1, Th. 4.2  and  the theorem on P. 535 of [W2], 
we immediately have the following:
\proclaim{Corollary 4.2}
Conj. 1 of 2.3 is true for the following inclusions:\par
(1) $G_{k_1+k_2+...+k_m} \subset G_{k_1} \times... \times
G_{k_m}$ where the inclusion is diagonal ,$k_i\in {\Bbb N}, i=1,...,m$
and $G=SU(n)$;\par
(2) $H_{lk} \subset G_l$, if  $H_{k} \subset G_1$
is a conformal inclusion where $k$
is the Dynkin index , $l\in {\Bbb N},l>1$,  $H$ and $G$ are simple and
of type $A$; \par
(3) $H\subset G_m$, where $H$ is the Cartan subgroup of $G$ , 
$m\in {\Bbb N},m>1$ and
G is a simple type A group.
\endproclaim   
Note in (2) and (3) of Cor. 4.2, we restrict $l>1$ and $m>1$ respectively
to avoid the trivial case of conformal inclusions. 
\subheading{4.3  
$SU(N) \subset SU(N)_{m'} \times SU(N)_{m''}$}
In this subsection we consider the coset $H\subset G_L$ with  
$H:=SU(N), G_L:=SU(N)_{m'} \times SU(N)_{m''}$, where
the embedding $H\subset G_L$ is diagonal. 
Let $\Lambda_1,...,\Lambda_{N-1}$ be the fundamental weights of $sl(N)$.
Let $k\in \Bbb N$. Recall that the set of integrable weights of the affine
algebra $\widehat{sl(N)}$ at level $k$ is the following subset of the weight
lattice of  $sl(N)$:
$$
P_{++}^{(h)} = \{ \lambda = \lambda_1 \Lambda_1 +...+ \lambda_{N-1}
 \Lambda_{N-1} | \lambda_i\in \Bbb N,  \lambda_1+...+ \lambda_{N-1} <h
\}
$$ where $h=k+N$.  This set admits a $\Bbb {Z}_N$ automorphism generated  
by 
$$
\sigma_1: \lambda=( \lambda_1, \lambda_2,..., \lambda_{N-1}) \rightarrow
\sigma_1( \lambda) = 
(h-\sum_{j=1}^{N-1} \lambda_j, \lambda_1,...,\lambda_{N-2})
.$$  We define the color $\tau(\lambda):\equiv \sum_i (\lambda_i-1)i 
 \text{\rm mod} (N)$ and $Q$ to be the root lattice of $\widehat{sl(N)}$
(cf. \S1.3 of [KW]). Note that $\lambda \in Q$ if and only if 
$\tau(\lambda) \equiv 0 \ \text{\rm mod}(N)$.
\par 
As in 2.1, we 
use $i$ (resp. $\alpha$) to denote the irreducible positive
energy representations of $LG$ (resp. $LH$).  To
compare our notations with that of \S2.7 in  
[KW], note that our $i$ is $(\Lambda',\Lambda'')$
of [KW] , and our $\alpha$ is $\Lambda$ of [KW]. 
We will  identify $i=(\Lambda',\Lambda'')$ and $\alpha =\Lambda$
where $\Lambda',\Lambda''$, $\Lambda$ are the weights of 
$sl(N)$ at levels $m', m'', m'+m''$ respectively.
Denote by $0',0'',0$ the vacuum representations of 
$\widehat{sl(N)}$ at level $m',m''$ and $m'+m''$  respectively.
Note that
for $i= (\Lambda',\Lambda''), \alpha =\Lambda$,   
by Th. 1.2 of [X1] the statistical dimensions of $i$ and $\alpha$
are given by  
$$
d_i =\frac{a(\Lambda') a(\Lambda'')}{a(0') a(0'')},
d_\alpha = \frac{a(\Lambda)}{a(0)}
$$
where the positive numbers $a(\Lambda),a(\Lambda')$ and $a(\Lambda'') $ 
are defined as in (0.4b) of [KW] ($a(\Lambda)$ is also equal to 
$S_{\Lambda_0}^{(\Lambda)}$ as defined on P. 362 of [X1]).
\par

Suppose 
$$
i=({\Lambda_1}',{\Lambda_1}''), j=({\Lambda_2}',{\Lambda_2}''),
k=({\Lambda_3}',{\Lambda_3}''), \alpha= {\Lambda_1}, \beta={\Lambda_2},
\delta=\Lambda_3.
$$ 
Then the fusion coefficients
$N_{ij}^k := N_{{\Lambda_1}'{\Lambda_2}'}^{\Lambda_3'}
 N_{{\Lambda_1}''{\Lambda_2}''}^{{\Lambda_3}''}$  
(resp. $N_{\alpha\beta}^\delta := N_{\Lambda_1\Lambda_2}^{\Lambda_3} $ 
)of $LG$ (resp. $LH$)
are well known and they are given by Verlinde formula (cf. Cor. 1 on 
P. 536 of [W2] and P. 288 of [Kac]).
\par
Recall $\pi_{i,\alpha}$ are the covariant representations of
the coset $H\subset G_L$. 
The set of all $(i,\alpha):= (\Lambda',\Lambda'',\Lambda)$ 
which appear 
in the decompositions of $\pi^i$ of $LG$ with respect to $LH$
is denoted by $exp$. This set is determined on P. 194 of
[KW] to be $(\Lambda',\Lambda'',\Lambda) \in exp$ iff 
$\Lambda'+\Lambda''-\Lambda \in Q$. 
The  ${\Bbb Z_N}$ action on $(i,\alpha), \forall i,
\forall \alpha$  is denoted by 
$$
\sigma(i,\alpha):= (\sigma(i), \sigma(\alpha))
=(\sigma (\Lambda'), \sigma (\Lambda''), \sigma (\Lambda))
, \sigma \in  {\Bbb Z_N}. $$
This is also known as diagram automorphisms
since they corresponds to the automorphisms of Dynkin diagrams.  
Note that $d_{\sigma (i)} = d_i, d_{\sigma (\alpha)} = d_\alpha$
by (3.2) of [Wal] and the formula for statistical dimensions above.
Also note that
this  ${\Bbb Z_N}$ action preserves $exp$ and therefore induces a 
${\Bbb Z_N}$ action on  $exp$. For each $(i,\alpha) \in exp$, we will 
denote by $[i,\alpha]$ its orbit in $exp$ under the  ${\Bbb Z_N}$ action.

\proclaim{Theorem 4.3}
Let $H\subset G_L$ be as in the previous paragraph. Then:\par 
(1)
$a_{(i,\alpha)\otimes 1} = \sigma_i a_{1\otimes \bar{\alpha}}$
for all $(i,\alpha) \in exp$; \par
(2) Assume the action of $Z_N$ on  $exp$ is faithful, i.e., if
$\sigma(i)=i, \sigma(\alpha) = \alpha$ for some 
$(i,\alpha)\in exp$, then $\sigma=id$.
Then  the covariant 
representations $\pi_{i,\alpha}$ are irreducible and 
$\pi_{i,\alpha}$ is unitarily equivalent to $\pi_{j,\beta}$ as 
covariant representations iff
$\sigma(i)=j,\sigma(\alpha) = \beta$ for some $\sigma\in {\Bbb Z_N}$;\par    
(3) Suppose the conditions of (2) hold.
Denote by $EXP$ the set of all irreducible localized sectors 
corresponding to the the covariant 
representations $\pi_{i,\alpha}, \forall (i,\alpha)\in exp$.
Then the set $EXP$  is  in one to one
correspondence with the set $\{ [i,\alpha],\forall (i,\alpha)\in exp \}$. 
Denote the elements of $EXP$  by $[i,\alpha]$. 
Define 
$$
C_{[i,\alpha][j,\beta]}^{[k,\delta]} :
= \sum_{\sigma \in {\Bbb Z_N}} N_{ij}^{\sigma(k)}
N_{\alpha \beta}^{\sigma(\delta)}
.$$ Then the compositions  are given by:
$$
[i,\alpha][j,\beta] = \sum_{[k,\delta]} C_{[i,\alpha][j,\beta]}^{[k,\delta]}
[k,\delta]
;$$ \par
(4) Conj. 2 of 2.3 is true for $H\subset G_L$.
\endproclaim
\demo{Proof}
By Cor. 4.2 and Th. 4.1, for any $(i,\alpha), (i',\alpha')$ we have:
$$
\align
\langle a_{1\otimes \bar{\alpha}}\sigma_i,  a_{1\otimes \overline{\alpha'} }
\sigma_{i'}
\rangle &=  \langle a_{1\otimes \bar{\alpha}}a_{1\otimes {\alpha'}}
,\sigma_{\bar i}\sigma_{i'} \rangle \\
&=\langle a_{1\otimes \bar{\alpha} {\alpha'}}
,\sigma_{\bar i}\sigma_{i'} \rangle \\
&=\langle \sum_\beta N_{\bar{\alpha} \alpha'}^\beta a_{1\otimes \beta}
,\sum_j N_{\bar{i} i'}^j \sigma_j \rangle \\
&=\langle \sum_\beta N_{\bar{\alpha} \alpha'}^\beta 1\otimes \beta
,\sum_{j,\delta}N_{\bar {i} i'}^j (j,\delta) \otimes \delta \rangle \\
&= \sum_\beta N_{\bar{\alpha} \alpha'}^\beta \sum_jN_{\bar i i'}^j \langle   1
,(j,\beta)  \rangle, \\
\endalign
$$ where $1$ in the last = stands for the identity sector of
$A(I)$, which by Th. 2.3, corresponds to the representation
$\pi_{0,0}$. So
$$
 \langle   1
, (j,\beta)  \rangle \neq 0
$$
if and only if 
$\pi_{0,0}$ appears as an irreducible summand in  
$\pi_{j,\beta}$ by the remarks after lemma 4.2.
%Recall that
%$$
%H^j=\sum_{\alpha} H_{j,\alpha}\otimes H_\alpha
%$$ 
Note that the 
vacuum vector (unique up to  a nonzero scalar)
 of $\pi_{0,0}$ has lowest energy (the eigenvalue of
the generator of the rotation group) $0$, and 
$\pi_{0,0}$ is the unique (up to 
unitary equivalence) irreducible representation  
with this property (cf. remarks after Prop. 2.1). So
$\pi_{0,0}$ appears  as an irreducible summand in  
$\pi_{j,\beta}$ if and only if there exists
a nonzero (vacuum) vector in $H_{j,\beta}$ with lowest energy $0$.  The
set of such $(j,\beta)$ was introduced  on  P. 186 
 in [KW] with our 
$(j,\beta)$ corresponds to $(M,\mu \text{\rm mod} (\delta))$ in 
the notation of [KW] (in
the notation on P. 186 of [KW], $h_M- h_\mu\geq 0$ is the
eigenvalue of the generator of the rotation group in the coset 
Hilbert space by definition). 
This set in our case of diagonal inclusions
is determined in   
(2.7.12) of [KW]. Translate (2.7.12) of [KW] into the
notations of this paper,   the statement is that $\pi_{0,0}$ 
appears as an irreducible summand in  
$\pi_{j,\beta}$ 
iff  there exists $\sigma \in {\Bbb Z_N}$ such
that $\sigma (0,0) = (j,\beta)$, where $0$ is used to denote the
vacuum representation of $G$ and $H$. 
Since  (cf. (3.3) of [Wal]) 
$$
N_{\bar i i'}^{\sigma(0)} = \delta_{\sigma (i),i'},
N_{\bar{\alpha} \alpha'}^{\sigma(0)} = \delta_{\sigma (\alpha),\alpha'}
,$$ we have:
$$
\langle a_{1\otimes \bar{\alpha}}\sigma_i,  a_{1\otimes \overline{\alpha'}}
\sigma_{i'}
\rangle = \sum_{\sigma \in {\Bbb Z}_N} 
\delta_{\sigma (i),i'} \delta_{\sigma (\alpha),\alpha'}
\langle \sigma (0,0),(0,0)\rangle. \tag **
$$ 
Ad (1): We will prove (1) by using an ``exhaustion'' trick similar
to the one used in \S 3 of [X3]. \par
Suppose $(i,\alpha) =(i',\alpha') = (0,\alpha)$,
where $0$ denotes the vacuum sector of $LG$.  Note
 that $\sigma(0)=0$ iff $\sigma =1$. So we conclude 
from (**) that $ a_{1\otimes \bar{\alpha}}$ is irreducible. By (2)
of Prop. 4.2, $ a_{1\otimes \bar{\alpha}} = a_{(0,\alpha)\otimes 1}$
if $(0,\alpha) \in exp$, and so $d_{(0,\alpha)} = d_\alpha$. \par 
 Note that (cf. the paragraph before Prop. 4.2) for fixed $i$, the
statistical dimension of 
$\rho \sigma_i \bar{\rho}$ is given by
$\sum_{\alpha} d_{i,\alpha} d_\alpha$ where the sum is over
those $\alpha$ with $(i,\alpha)\in exp$, which will be denoted by
$exp_i$. Note that if $i=(\Lambda', \Lambda''), \alpha=\Lambda$,
then $exp_i$ is a 
congruence class of $P_{++}^{m'+m''+N} \text{\rm mod} Q$
(congruent to $\Lambda' +\Lambda'' $) by definition. 
So $d_\rho^2 = \frac{1}{d_i} \sum_{\alpha\in exp_i} d_{i,\alpha} d_\alpha$.
By (3) of Prop. 4.2, $d_{(i,\alpha)} \leq d_i d_\alpha,
\forall (i,\alpha)\in exp$, hence
$$
d_\rho^2 \leq  \sum_{\alpha\in exp_i} d_{\alpha} d_\alpha
= \sum_{\alpha\in exp_0}  d_{\alpha} d_\alpha 
,$$ where tha last = follows from Cor. 2.7 of [KW]\footnotemark\footnotetext{
Note that our $\alpha$ correspond to $\Lambda$ in Cor. 2.7 of [KW], 
$d_\alpha=\frac{a(\Lambda)}{a(0)}$ where $0$ denotes the vacuum
representation, and $exp_i$ is a 
congruence class of $P_{++}^{m'+m''+N} \text{\rm mod} Q$
.}. But since
$d_{(0,\alpha)} = d_\alpha$, we have:
$$
d_\rho^2 = \sum_{\alpha\in exp_0}  d_{\alpha} d_\alpha
.$$ It follows that all the $\leq$'s above are actually $=$, in particular
$d_{i,\alpha} = d_i d_\alpha$, and
it follows from (2) of Prop. 4.2 that 
$a_{(i,\alpha)\otimes 1} = \sigma_i a_{1\otimes \bar{\alpha}}$.\par
So we have (see the paragraph before Th. 4.3): 
$$
d_{\sigma (0,0)}= d_{\sigma (0)} d_{\sigma (0)}=1.
$$ On the
otherhand 
$$
\langle \sigma (0,0),(0,0)\rangle \geq 1,
$$
it follows by comparing statistical dimensions that
$$
\langle \sigma (0,0),(0,0)\rangle = 1.
$$
So we can improve (**) to 
$$
\langle a_{1\otimes \bar{\alpha}}\sigma_i,  a_{1\otimes \overline{\alpha'} }
\sigma_{i'}
\rangle = \sum_{\sigma \in {\Bbb Z}_N} 
\delta_{\sigma (i),i'} \delta_{\sigma (\alpha),\alpha'}
. \tag *
$$

Ad (2): By assumption expression (*) holds
, for a unique $\sigma$ since the action is faithful. If $(i,\alpha) =
(i',\alpha')$, then $\sigma(i)=i',\sigma(\alpha) =\alpha'$ 
iff $\sigma =id$, and we conclude from (*) that
$$
\langle a_{1\otimes \bar{\alpha}}\sigma_i,  a_{1\otimes \bar{\alpha}}
\sigma_{i}
\rangle =1
,$$ i.e., $ a_{1\otimes \bar{\alpha}}\sigma_i$ is irreducible. By 
(1) of Prop. 4.2 and (1) of Th. 4.3 
we conclude that if $(i,\alpha)\in exp$, then
$(i,\alpha)$ is irreducible, i.e., $\pi_{i,\alpha}$ is irreducible
by the remarks after lemma 4.2. \par
If  $\sigma(i)=j,\sigma(\alpha) =\beta$, then 
$$ 1\geq \langle a_{1\otimes \bar{\alpha}}\sigma_i,  
a_{1\otimes \bar{\beta} }
\sigma_{j}\rangle \geq 1
,$$ where the first $\geq$ follows from the fact that 
$ a_{1\otimes \bar{\alpha}}\sigma_i,  a_{1\otimes \bar{\beta} }
\sigma_{j}$ are irreducible, and 
the second $\geq$ follows from (*).  So we must have
$$
 a_{1\otimes \bar{\alpha}}\sigma_i = a_{1\otimes \bar{\beta} }
\sigma_{j}.$$  By (1) of Prop. 4.2 and (1) of Th. 4.3 
$
\langle (i,\alpha),  (j,\beta) \rangle =1,
$ so
$(i,\alpha)$ is identical
to $(j,\beta)$ as sectors since both sectors are
irreducible. Hence $\pi_{i,\alpha}$ is unitarily equivalent to $\pi_{j,\beta}$
as covariant representations by the remarks after lemma 4.2.  
On the other hand if 
$\pi_{i,\alpha}$ is unitarily equivalent to $\pi_{j,\beta}$
as covariant representations, then $(i,\alpha)$ is equal to
$(j,\beta)$ as sectors by the remarks after lemma 4.2.
By (1) of Prop. 4.2 and (1) of Th. 4.3
$$
a_{1\otimes \bar{\alpha}}\sigma_i = a_{1\otimes \bar{\beta} }
\sigma_{j},$$
so we must have 
$\sigma (i)=j,\sigma (\alpha) = \beta$ for some $\sigma \in {\Bbb Z}_N$ by 
(*). This proves (2).
\par
Ad (3): Note  that if $(i,\alpha) \in exp,
(j,\beta) \in exp$, by (1) and Th. 4.1 we have:
$$
\align
a_{(i,\alpha) \otimes 1} a_{(j,\beta) \otimes 1} &= 
\sigma_i \sigma_j a_{1\otimes \bar \alpha}a_{1\otimes \bar \beta} \\
&= \sum_{k,\delta} N_{ij}^k N_{\alpha \beta}^\delta \sigma_k 
a_{1\otimes \bar \delta}
\endalign
$$   Note by the characterization of $exp$ (cf. P. 194 of [KW]) and
4.3.4 of [FKW], if  $N_{ij}^k N_{\alpha \beta}^\delta \neq 0$,
then $(k,\delta)\in exp$. Using (1) again we have:
$$
a_{(i,\alpha) \otimes 1} a_{(j,\beta) \otimes 1}= 
\sum_{(k,\delta)\in exp} N_{ij}^k N_{\alpha \beta}^{\delta} 
a_{(k, \delta) \otimes 1} 
.$$  (3) now follows from (1) of Prop. 4.2 and (2). \par

Ad (4): By (1) we need to show that
$$
\frac{b(i,\alpha)}{b(0,0)} = d_i d_\alpha,
$$
where under the identification of $i$ with $(\Lambda', \Lambda'')$ and
$\alpha$ with $\Lambda$, $b(i,\alpha)$ is given by (2.7.14) of [KW],
and use the notations there we have
$$
b(i,\alpha)= b(\Lambda', \Lambda''; \Lambda) = N a(\Lambda')  a(\Lambda'')
a(\Lambda).
$$
Note that $d(\Lambda)= \frac{a(\Lambda)}{a(0)}$, and we have
similar identities when $\Lambda, 0$ is replaced by 
$\Lambda', 0'$ and $\Lambda'', 0''$.   The proof of (4)
is now complete by definitions.
\enddemo \hfill \qed
\par
Note  that by (1) of Th. 4.3 and formula (*) 
if $\sigma(i)=i,\sigma(\alpha) =\alpha$ for some $\sigma \neq id$, 
$(i,\alpha)\in exp$, then
$a_{(i,\alpha) \otimes 1}$ is not irreducible. Hence $(i,\alpha)$
is not irreducible by (1) of Prop. 4.2. 
This happens for an example when $N=2, m'\in 2 {\Bbb N}$ 
and $m''=2$, which is related to supersymmetric conformal algebras
(cf. P. 195 of [KW]). \par   
Recall from subsection 2.3 coset  $W_N$ algebras 
with critical parameters  are defined to be the irreducible conformal
precosheaves of cosets 
$$
SU(N)_{m+1} \subset SU(N)_m \times SU(N)_1,
$$ which obviously satisfies
(1) to (4) of Th. 4.3. 
For $N=2$, Th. 4.3 is obtained in [Luke] by different methods.
So we obtain a proof of a long-standing
conjecture about representations of  coset  $W_N$ ($N>2$) algebras
with critical parameters, which is similar in the type $A$ case to
the conjecture 3.4
and Th. 4.3 stated  in \S 3 and 4 of [FKW]. 
To compare (3) of Th. 4.3 to Th. 4.3 of [FKW], note that $(p,p')$ in
[FKW] is identified with $(m',m'+1)$ in our case,  
and $(\lambda_i,
\lambda_i')$ in Th. 4.3 of [FKW] is identified with 
$[\Lambda_i', {\Lambda_i}'', \Lambda_i]$ in our case with
$\lambda_i =\Lambda_i', \lambda_i' = \Lambda_i$, and ${\Lambda_i}''$
is the unique element such that $ \Lambda_i'+{\Lambda_i}''-\Lambda_i
\in Q$.  Then the compostions in Th. 4.3 of [FKW] can be easily
checked to be the same as (3) of Th. 4.3 under the assumptions of
Th. 4.3 in [FKW].\par
%We end this subsection by noting that for $G$ of types other than $A$, similar
%results as Th. 4.3 can be proved.  We will discuss this elsewhere.   
\subheading{4.4 More examples} 
Let us consider the case $H\subset G_m$ with $G=SU(l)$ and
$H$ is the Cartan subalgebra of $G$, a $l-1$ dimensional torus.
This coset does not verify the conditions of Th. 2.3, but Th. 2.3
applies to this coset by the remark before Cor. 3.1. So we can apply
the results of section 4. \par
Recall the  equivalent relation  on ${\Bbb Z}^{l-1}$  in 3.1 
defined by:
$$ 
(n_1,...,n_{l-1}) \sim  (n_1',...,n_{l-1}')
$$ iff
there exists $(m_1,...m_{l-1}) \in {\Bbb Z}^{l-1}$ with
$m_1+...+m_{l-1} \in l{\Bbb Z}$ such that
$$(n_1',...,n_{l-1}') = (n_1+mm_1,...n_{l-1}+mm_{l-1})
.$$ Denote
the equivalence class of $(n_1,...,n_{l-1})$ by $[n]:=[n_1,...,n_{l-1}]$
or simply $[n]$,
and they are   used to denote irreducible projective representations of
$LH$.\par
Let $\Lambda$ be an irreducible weight of
$SU(l)$ at level $m$. Let $\tau(\Lambda)$ be the 
color of $\Lambda$ , see 4.3. Then:
$$
H^\Lambda = \oplus_{ \{[n] |(\sum_in_i - \tau (\Lambda)) \in l{\Bbb Z}
\}} H_{\Lambda, [n]} \otimes H_{[n]}
,$$ cf. \S2.6 of [KW]. So the set $exp= \{ (\Lambda, [n]) ,
|(\sum_in_i - \tau (\Lambda))\in l{\Bbb Z} \} $ is determined, and 
$(i,\alpha)= (\Lambda, [n])$ in 
the present case.  Note that 
the sector $[n]$ corresponds to an automorphism and has statistical dimension
$1$ (cf. section 3 before Cor. 3.1). 
So
$a_{1\otimes 
[n]}$ has  statistical dimension
$1$, and it follows that
$ \overline{a_{1\otimes [n]}} a_{1\otimes [n]}$ is the identity sector. So   
$$
\langle a_{1\otimes 
[n]} \sigma_\Lambda, a_{1\otimes 
[n]} \sigma_\Lambda \rangle 
= \langle \overline{a_{1\otimes [n]}}a_{1\otimes [n]}\sigma_\Lambda, 
\sigma_\Lambda \rangle =  \langle \sigma_\Lambda, \sigma_\Lambda \rangle
=1.
$$
Hence  $ a_{1\otimes 
[n]} \sigma_\Lambda   
$ is irreducible for any $[n], \Lambda$.
So by Cor. 3.1 and (2) of Prop. 4.2 if $(\Lambda, [n]) \in exp$,
$$
a_{(\Lambda, [n])\otimes 1} = \overline{a_{1\otimes [n]}} \sigma_\Lambda,
$$ and by (1) of Prop. 4.2 we have the following fusion rules:
$$
(\Lambda, [n]) (\Lambda', [n']) = \sum_{\Lambda''}
N_{\Lambda \Lambda'}^{\Lambda''}
(\Lambda'', [n+n'])
,$$ for any $(\Lambda, [n]), (\Lambda', [n']) \in exp$.
By using (2.6.14) of [KW], one can immediately check that 
Conj. 2 of 2.3 is  true in this case, as we did in the proof of 
(4) of Th. 4.3.  
Note this coset is related to Parafermions, see [DL] for
an approach using vertex operator algebras. \par
Next we consider the case $SU(2)_{4k} \subset SU(3)_k$. 
When $k=1$, 
$SU(2)_4 \subset SU(3)_1$ is a conformal inclusion, so by
Cor. 4.2, the coset $SU(2)_{4k}\subset SU(3)_k, k>1$ is rational, i.e.,
it verifies Conj. 1 of  2.3. In this particular example, we
will use the notations in \S6 of [DJ1]. Thus the weights of
$\widehat{sl(2)}$ are labeled by an integer $l$ with $0$ 
denoting the vacuum representation and the weights of 
$\widehat{sl(3)}$ are labeled by a pair of integers $(pq)$ with
$(00)$ denoting the vacuum representation. So $(i,\alpha):= (pq,l)$. \par
%We start with the case $k=3$. In this case $((pq),l) \in exp$ iff
%$l$ is even (cf. P.4115 of [DJ1]). By calculating the 
%lowest weights it is easy to show that:
%$$
%\langle (00,0), (pq,l) \rangle =0
%$$ unless $(pq,l)= (00,12)$. By using Th. 4.1 and Prop.4.2 we 
%find 
%$$
%a_{(pq,l)\otimes 1} = \sigma_{pq} a_{1\otimes l}, \forall
%(pq,l) \in exp
%,$$ 
%and so the compositions of $(pq,l)$'s are completely
%determined by (1) of Prop.4.2. Also $(pq,l)=(pq,12-l)$ can
%be checked almost immediately which corresponds to
%the field identification (6.2) of [DJ1]. \par
In general, the
identifications of $(i,\alpha) =(j,\beta)$ as sectors are  related 
to certain Dynkin diagram automorphisms,  as in  (2) of
Th. 4.3. See [SY] for more examples.  However, in [DJ1],
a ``Maverick" coset is given which violates the above identification rule.
This is the coset $SU(2)_8 \subset SU(3)_2$. From the first line of 
table on  P. 4117 of [DJ1], we have:
$$
(00,0)=(00,8)=(11,4) 
$$ as sectors by Th. 2.3 
and the  fact that
the vacuum representation is the unique (up to unitary equivalence)
representation which contains  a nonzero (vacuum) vector with 
lowest energy $0$. 
By using the above equation, Th. 4.1 and 
Prop. 4.2
we find that the irreducible subsectors of $a_{(pq,l)\otimes 1}$
generate a $6$ dimensional ring whose basis is:
$$
\align
1, x:=a_{(00,4)\otimes 1} &= a_{1\otimes 4} -\sigma_{11},\\
y:=a_{(10,2)} & = \sigma_{10} a_{1\otimes 2} - \sigma_{02} a_{1\otimes 2}
,\\
\bar y:=a_{(01,2)} & = \overline{a_{(10,2)}}, \\
z:=a_{(10,4)} & = \sigma_{02},
\bar z:=a_{(01,4)}  = \overline{a_{(10,4)}}
\endalign
$$  where $\bar a$ is the conjugate of $a$.  The statistical dimensions 
of  $ a_{(00,4)\otimes 1},  a_{(10,2)}, 
a_{(10,4)}$ are given by $\frac{\sqrt 5+1}{2}, \frac{\sqrt 5+1}{2},1$
respectively.
The fusion rules are also completely determined by the above formula,
and we have 
$$
[x^2]=[1]+[x], [y\bar y]=[1]+[x], [z^3]=[1], [y]= [xz]
.$$ 
By (1) of Prop. 4.2, the above formula determines the 
structure of the ring generated by the coset sectors.\par 
Conj. 2. of 2.3 is also  easily verified in this case.
For more such  ``Maverick" cosets, see [DJ2]. \par
\heading References \endheading
\roster
\item"{[ABI]}" D. Altschuler, M. Bauer and C. Itzykson, {\it The
branching rules of conformal embeddings},  Comm. Math. Phys.,
 {\bf 132}, 349-364
(1990).  
%\item"{[FKW]}" E. Frenkel, V. Kac and M. Wakimoto, {\it Characters and
%Fusion rules for W-algebras via Quantized Drinfeld-Sokolov Reductions},
%RIMS-861, 1992.
\item"{[B]}" R. Borcherds, {\it Vertex algebras, Kac-Moody
algebras and the Monster,}, Proc. Nat. Acad. Sci., U.S.A. 83, 3068-3071
(1986).
%\item"{[BB]}" F. Alexander Bais and P. Bouwknegt, {\it 
%A classification of subgroup truncations of the bosonic string},
%Nucl.Phys. B 279, 561-70 (1987).
\item"{[BBSS]}" F.A. Bais, P. Bouwknegt, K. Schoutens, M. Surridge,
{\it Coset construction for extended Virasoro algebras,} Nucl. Phys. B
304, 371-391, (1988).
%\item"{[BE1]}" J. B\"{o}ckenhauer and D. Evans, {\it Modular
%Invariants, Graphs and $\alpha$-Induction for Nets of Subfactors I,
%} Comm. Math. Phys., {\bf 197}, 361-386 (1997).
%\item"{[BE2]}" J. B\"{o}ckenhauer and D. Evans, {\it Modular
%Invariants, Graphs and $\alpha$-Induction for Nets of Subfactors II,
%} Comm. Math. Phys., {\bf 200}, 57-103 (1999).
\item"{[D]}" C. Dong, {\it Introduction to vertex operator algebras I},
\par
$S\bar{u}rikaisekikenky\bar{u}sho \ \ K\bar{o}ky\bar{u}roku$, 
No. 904 (1995), 1-25. Also
see q-alg/9504017.
\item"{[DL]}" C. Dong and J. Lepowsky, {\it Generalized vertex algebras
and relative vertex operators,} Progress in Mathematics, 112 (1993).
\item"[DJ1]" D. Dunbar and K. Joshi,
{\it Characters for Coset conformal field theories and Maverick
examples}, Inter. J. Mod. Phys. A, Vol.8, No. 23 (1993), 4103-4121.
\item"[DJ2]"  D. Dunbar and K. Joshi, {\it Maverick examples of 
Coset conformal field theories}, 
Mod. Phys. Letters A, Vol.8, No. 29 (1993), 2803-2814.
\item"{[Dix]}" J. Dixmier, {\it von Neumann Algebras}, North-Holland
Publishing Company, 1981.
\item"{[FG]}" J. Fr\"{o}hlich and F. Gabbiani, {\it Operator algebras and
Conformal field theory, } Comm. Math. Phys., {\bf 155}, 569-640 (1993).
\item"{[FKW]}" E. Frenkel, V. Kac and M. Wakimoto, {\it Characters and
Fusion rules for W-algebras via Quantized Drinfeld-Sokolov Reductions,}
Comm. Math. Phys., {\bf 147}, 295-328 (1992).
\item"{[FLM]}" I. B. Frenkel, J. Lepowsky and J. Ries,
{\it Vertex operator algebras and the Monster}, Academic, New York, 1988.
\item"{[FZ]}" I. Frenkel and Y. Zhu,
{\it Vertex operator algebras associated to representations of
affine and Virasoro algebras,} Duke Math. Journal (1992), Vol. 66, No. 1,
123-168.
\item"{[GKO]}" P. Goddard, A. Kent  and D. Olive,  { \it 
Unitary representations of 
Virasoro and super-Virasoro algebras,} 
Comm. Math. Phys. 103 (1986), No. 1, 105-119.
\item"{[GL]}"  D. Guido and R. Longo, {\it  The conformal  spin and
statistics theorem},  \par
Comm.Math.Phys., {\bf 181}, 11-35 (1996).  
%\item"{[GL2]}"  D. Guido and R.Longo, {\it  Relativistic invariance and charge
%conjugation in quantum field theory},  \par
%Comm.Math.Phys., {\bf 148}, 521-551 (1992)  
\item"[GW]" R. Goodman and N. Wallach, {\it Structure and unitary cocycle
representations of loop groups and the group of diffeomorphisms of the
circle}, J. Reine Angew. Math 347 (1984) 69-133.
\item"{[Haag]}" R. Haag, {\it Local Quantum Physics}, Springer-Verlag 1992.
\item"{[J]}" V. Jones, {\it Fusion en al\'gebres de Von Neumann et groupes
de lacets (d'apr\'es A. Wassermann),} Seminarie Bourbaki, 800, 1-20,1995.
\item"{[KW]}"  V. G. Kac and M. Wakimoto, {\it Modular and conformal
invariance constraints in representation theory of affine algebras},  
Advances in Math., {\bf 70}, 156-234 (1988).
\item"{[Kacv]}"  V. G. Kac, {\it Vertex algebras for beginners}, AMS,
1997.
\item"{[Kac]}"  V. G. Kac, {\it Infinite dimensional Lie algebras}, 
3rd Edition,
Cambridge University Press, 1990. 
\item"{[KS]}" D. Karabali and H. J. Schnitzer, {\it BRST quantization
of the gauged WZW action and coset conformal field theories,}
Nucl. Phys. B 329 (1990), no. 3, 649-666.
\item"[KT]"  A. Tsuchiya and Y. Kanie, {\it Vertex Operators in conformal
field theory on $P^1$ and monodromy representations of braid group,}
Adv. Studies in Pure Math. 16 (88), 297-372.
\item"{[L1]}"  R. Longo, Proceedings of International Congress of
Mathematicians, 1281-1291 (1994).
\item"{[L2]}"  R. Longo, {\it Duality for Hopf algebras and for
subfactors},
I, Comm. Math. Phys., {\bf 159}, 133-150 (1994).
\item"{[L3]}"  R. Longo, {\it Index of subfactors and statistics of
quantum fields}, I, Comm. Math. Phys., {\bf 126}, 217-247 (1989).
\item"{[L4]}"  R. Longo, {\it Index of subfactors and statistics of
quantum fields}, II, Comm. Math. Phys., {\bf 130}, 285-309 (1990).
\item"{[L5]}"  R. Longo, {\it An analogue of the Kac-Wakimoto formula and
black hole conditional entropy,} 
Comm.Math.Phys., {\bf 186}, 451-479 (1997).
\item"{[L6]}"  R. Longo, {\it Minimal index and braided subfactors,
} J. Funct.Analysis {\bf 109} (1992), 98-112.
\item"{[LR]}"  R. Longo and K.-H. Rehren, {\it Nets of subfactors},
Rev. Math. Phys., {\bf 7}, 567-597 (1995).  
\item"{[Luke]}" T. Luke, {\it Operator algebras and Conformal Field Theory
of the Discrete series representations of Diff($S^1$),} Dissertation,
Cambridge (1994).
\item"[MT]"  M. Takesaki, {\it Conditional expectation in von Neumann
algebra,} J. Funct. Analysis 9 (1972), 306-321.
\item"{[MS]}" G. Moore and N. Seiberg, {\it Taming the conformal zoo},
Lett. Phys. B , {\bf 220}, 422-430, (1989).
\item"[Mv]" F. J. Murray and J. v. Neumann, {\it On rings of Operators},
Ann. Math. 37 (1936), 116-229.
\item"[Po]"  S. Popa, {\it Classification of amenable subfactors of type
II}, Acta Math.172 (1994), 352-445.
\item"{[PP]}" M. Pimsner and S. Popa, {\it Entropy and index for subfactors,} 
\par
Ann. \'{E}c.Norm.Sup. {\bf 19}, 57-106 (1986). 
\item"[PS]" A. Pressly and G. Segal, {\it Loop Groups,} O.U.P. 1986.
\item"[RS]" M. Reed and B. Simon, 
{\it Methods of Mathematical Physics I: Functional Analysis},
Academic Press 1980.
\item"[Stra]" S. Stratila, {\it Modular Theory in Operator Algebras,}
Editura Academiei, 1981.
\item"[SY]" A. N. Schellekens and S. Yankielowicz, 
{\it Field identification
fixed points in the coset construction,}
Nucl. Phys. B 334, 67-102, (1990).
\item"[Wal]" M. Walton, {\it  Fusion rules of Wess-Zumino-Witten Models},
 Nuclear physics B, 340: 777-789, 1990.    
\item"[We]" H.Wenzl, {\it Hecke algebras of type A and subfactors},
Invent. Math. 92 (1988), 345-383.
\item"{[W1]}"  A. Wassermann, Proceedings of International Congress of
Mathematicians, 966-979 (1994).
\item"{[W2]}"  A. Wassermann, {\it Operator algebras and Conformal
field theories III},  Invent. Math. 133 (1998), 467-538.
\item"[Watts]" G.M.T. Watts, {\it W-algebras and coset models}, Physics
Lett.B 245, 65-71, 1990.
\item"{[Witten]}" E. Witten, {\it The central charge in three dimensions,}
530-559, V.G. Knizhnik's Memorial Volume, 1989.
\item"[X1]" F. Xu, {\it   New braided endomorphisms from conformal
inclusions, } \par
Comm.Math.Phys. 192 (1998) 349-403.
\item"[X2]" F. Xu, {\it Applications of braided endomorphisms from
conformal inclusions,} 
Inter. Math. Res. Notice., No.1, 5-23 (1998), see also q-alg/9708013, 
and Erratum, Inter. Math. Res. Notice., No.8, (1998). 
\item"[X3]" F. Xu, {\it Jones-Wassermann subfactors for 
disconnected intervals}, \par 
 q-alg/9704003.
\item"{[Y]}"  S. Yamagami, {\it A note on Ocneanu's approach to Jones
index theory}, Internat. J. Math., {\bf 4}, 859-871 (1993). 
\endroster 
\enddocument